\documentclass[reqno]{amsart}

\usepackage{geometry}
\usepackage{graphicx}
\usepackage{float}
\usepackage{hyperref}
\usepackage{ellipsis}
\usepackage{amsmath,amsthm}
\usepackage{enumerate}
\usepackage{comment}
\usepackage{subcaption}
\usepackage{url}
\usepackage{cite}

\allowdisplaybreaks

\graphicspath{{./images/}}

\newcommand{\Q}{\mathbb{Q}}
\newcommand{\R}{\mathbb{R}}
\newcommand{\C}{\mathbb{C}}
\newcommand{\T}{\mathbb{T}}
\newcommand{\Z}{\mathbb{Z}}
\newcommand{\inner}[1]{\langle #1 \rangle}
\newcommand{\floor}[1]{\lfloor #1 \rfloor}
\renewcommand{\vec}[1]{\mathbf{#1}}
\renewcommand{\(}{\left(}
\renewcommand{\)}{\right)}
\renewcommand{\phi}{\varphi}

\newtheorem{Proposition}{Proposition}
\newtheorem{Corollary}{Corollary}
\newtheorem{Theorem}{Theorem}
\newtheorem{Lemma}{Lemma}

\begin{document}

\title[Gauss' hidden menagerie: from cyclotomy to supercharacters]{Gauss' hidden menagerie:\\from cyclotomy to supercharacters}

\author[S.~R.~Garcia]{Stephan Ramon Garcia}
\address{\textsc{Department of Mathematics\\
	Pomona College\\
	Claremont, California\\
	91711\\
	USA}}
\email{Stephan.Garcia@pomona.edu}
\urladdr{\url{http://pages.pomona.edu/~sg064747}}
\thanks{Partially supported by National Science Foundation
	Grant DMS-1265973.}

\author[T.~Hyde]{Trevor Hyde}
\address{\textsc{Department of Mathematics\\
	University of Michigan\\
	2074 East Hall\\
	530 Church Street\\
	Ann Arbor, MI 48109\\
	USA}}
\email{tghyde@umich.edu}

\author[B.~Lutz]{Bob Lutz}
\address{\textsc{Department of Mathematics\\
	University of Michigan\\
	2074 East Hall\\
	530 Church Street\\
	Ann Arbor, MI  48109\\
	USA}}
\email{boblutz@umich.edu}
\urladdr{\url{http://www-personal.umich.edu/~boblutz}}


\maketitle

\vspace{-3mm}
\section{Introduction}

At the age of eighteen, Gauss established the constructibility of the 17-gon, a result that had eluded mathematicians for two millennia. At the heart of his argument was a keen study of certain sums of complex exponentials, known now as~\emph{Gaussian periods}. These sums play starring roles in applications both classical and modern, including Kummer's development of arithmetic in the cyclotomic integers~\cite{Kolmogorov} and the optimized AKS primality test of H.~W.~Lenstra and C.~Pomerance~\cite{Agrawal,Lenstra}. In a poetic twist, this recent application of Gaussian periods realizes ``Gauss' dream'' of an efficient algorithm for distinguishing prime numbers from composites~\cite{Granville}.

We seek here to study Gaussian periods from a graphical perspective. It turns out that these classical objects, when viewed appropriately, exhibit a dazzling and eclectic host of visual qualities. Some images contain discretized versions of familiar shapes, while others resemble natural phenomena. Many can be colorized to isolate certain features; for details, see Section~\ref{sec:cycsup}.

\begin{figure}[ht]
	\begin{subfigure}[b]{0.41\textwidth}
		\includegraphics[width=\textwidth]{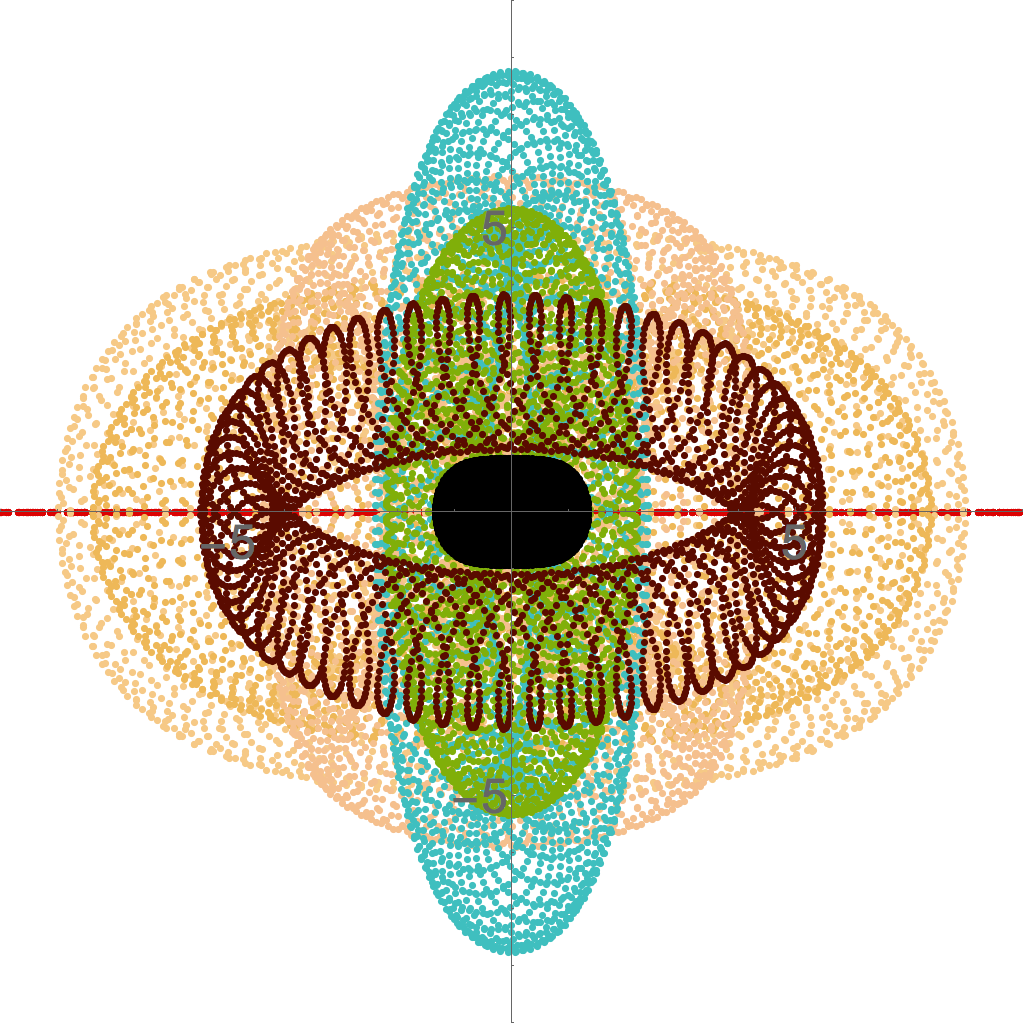}
		\caption{$n=29\cdot 109\cdot 113$, $\omega=8862$, $c=113$}
	\end{subfigure}
	\quad
	\begin{subfigure}[b]{0.41\textwidth}
		\includegraphics[width=\textwidth]{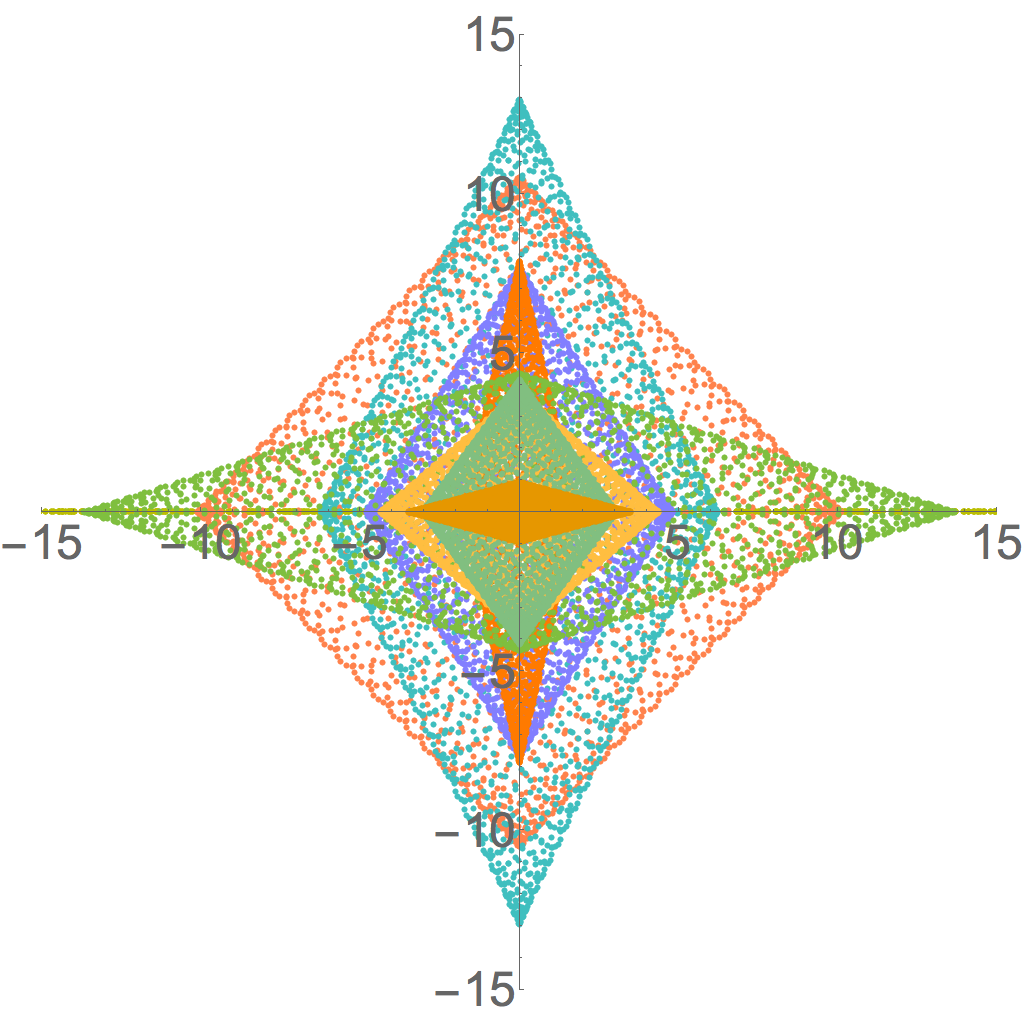}
		\caption{$n=37 \cdot 97 \cdot 113$, $\omega=5507$, $c=113$}
	\end{subfigure}
	\caption{Eye and jewel --- Images of \emph{cyclic supercharacters} correspond to sets of Gaussian periods. For notation and terminology, see Section~\ref{sec:cycsup}.}
	\label{fig:eye}
\end{figure}


\section{Historical context}

The problem of constructing a regular polygon with compass and straight-edge dates back to ancient times. Descartes and others knew that with only these tools on hand, the motivated geometer could draw, in principle, any segment whose length could be written as a finite composition of sums, products and square roots of rational numbers~\cite{Dunham}. Gauss' construction of the 17-gon relied on showing that
\[
	16\cos\(\frac{2\pi}{17}\)=
	-1+\sqrt{17}+\sqrt{34-2\sqrt{17}} +2\sqrt{17+3\sqrt{17}
	-\sqrt{34-2\sqrt{17}}-2\sqrt{34+2\sqrt{17}}}
\]
was such a length. After reducing the constructibility of the $n$-gon to drawing the length $\cos\(\frac{2\pi}{n}\)$, his result followed easily. So proud was Gauss of this discovery that he wrote about it throughout his career, purportedly requesting a 17-gon in place of his epitaph\footnote{H.~Weber makes a footnote of this anecdote in~\cite[p.~362]{Weber} but omits it curiously from later editions.}. While this request went unfulfilled, sculptor Fritz Schaper did include a 17-pointed star at the base of a monument to Gauss in Brunswick, where the latter was born~\cite{Lenstra2}.

\begin{figure}[ht]
	\begin{subfigure}[b]{0.484\textwidth}
		\includegraphics[width=\textwidth]{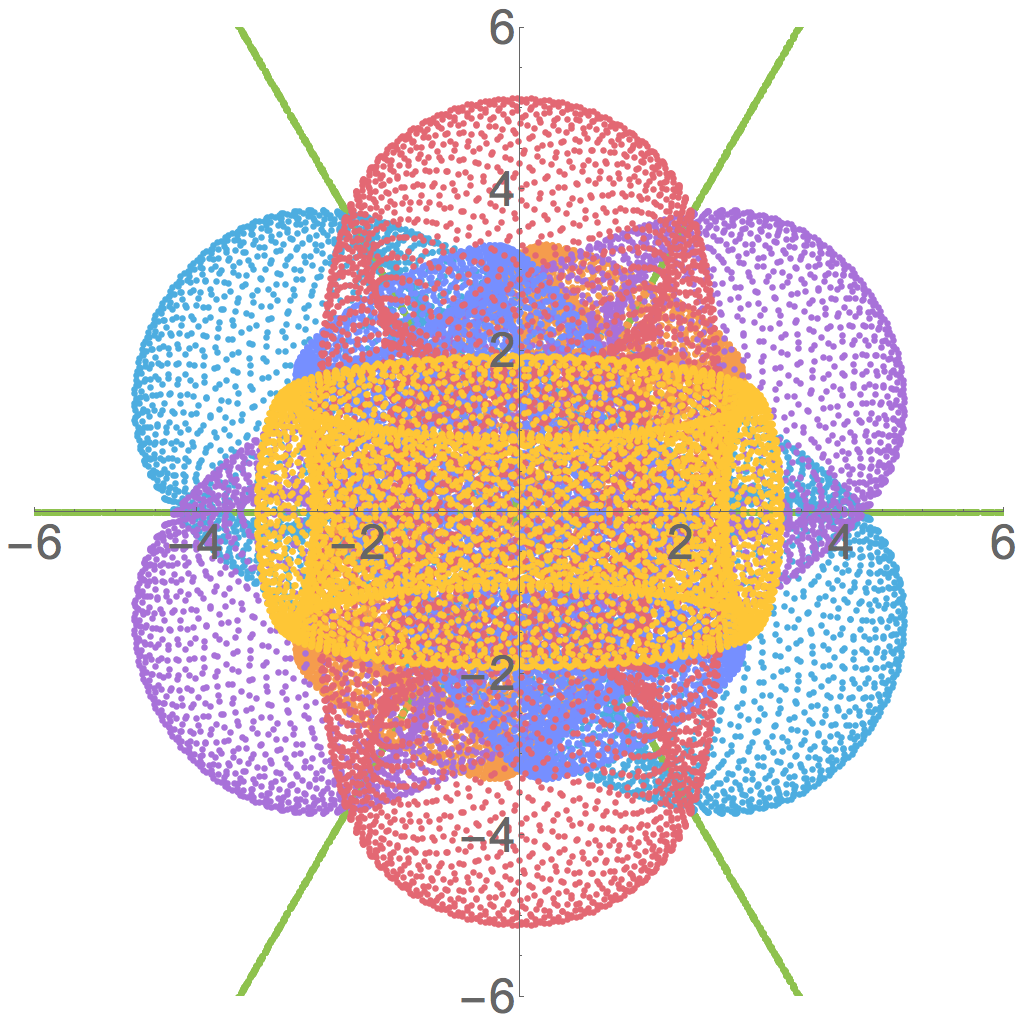}
		\caption{\scriptsize $n=3\cdot 5\cdot 17\cdot 29\cdot 37$, $\omega=184747$, $c=3\cdot 17$}
	\end{subfigure}
	\quad
	\begin{subfigure}[b]{0.484\textwidth}
		\includegraphics[width=\textwidth]{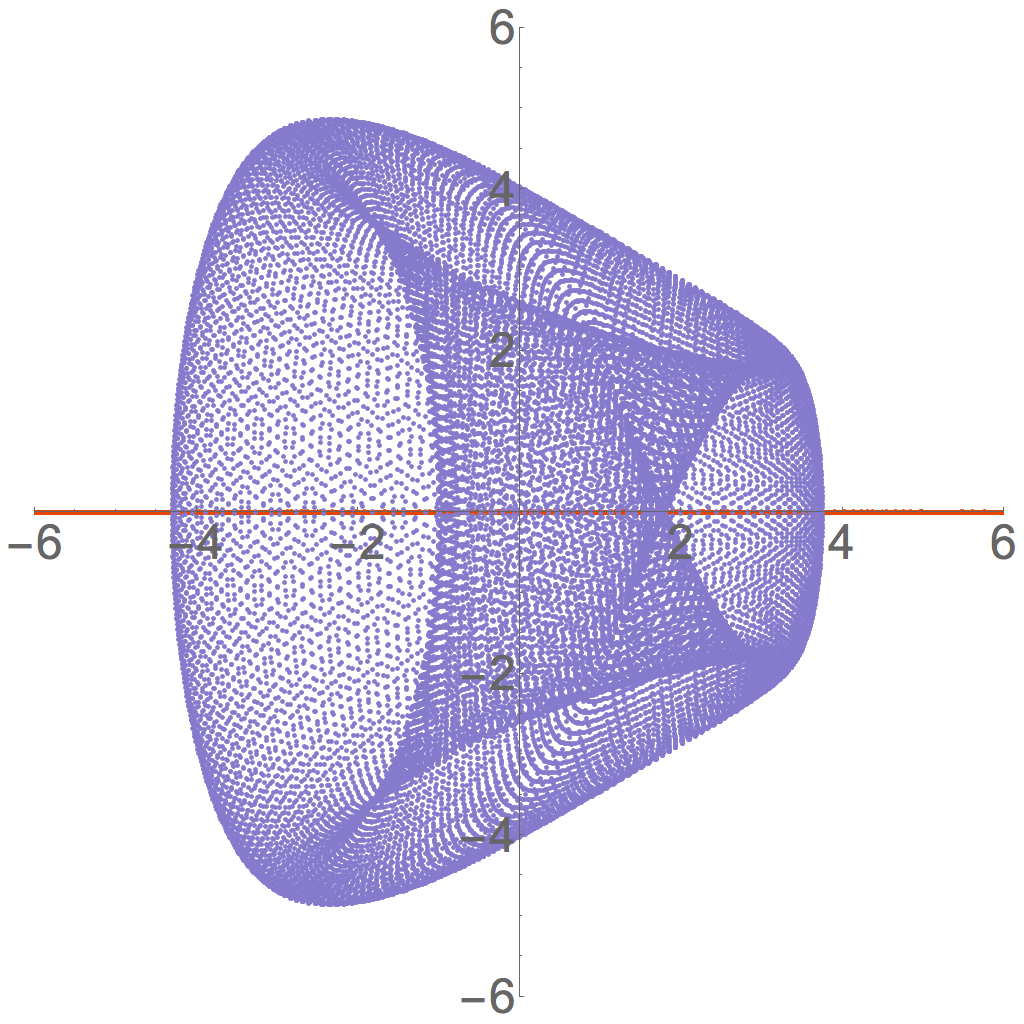}
		\caption{\scriptsize $n=13\cdot 127\cdot 199$, $\omega=6077$, $c=13$}
	\end{subfigure}
	\caption{Disco ball and loudspeaker --- Images of \emph{cyclic supercharacters} correspond to sets of Gaussian periods. For notation and terminology, see Section~\ref{sec:cycsup}.}
	\label{fig:crest}
\end{figure}

Gauss went on to demonstrate that a regular $n$-gon is constructible if Euler's totient $\phi(n)$ is a power of 2. He stopped short of proving that these are the only cases of constructibility; this remained unsettled until J.~Petersen completed a largely neglected argument of P.~Wantzel, nearly three quarters of a century later~\cite{Lutzen}. Nonetheless, the chapter containing Gauss' proof has persisted deservedly as perhaps the most well-known section of his \emph{Disquisitiones Arithmeticae}. Without the language of abstract algebra, Gauss initiated the study of \emph{cyclotomy}, literally ``circle cutting,'' from an algebraic point of view.

\begin{figure}[ht]
	\begin{subfigure}[b]{0.484\textwidth}
		\includegraphics[width=\textwidth]{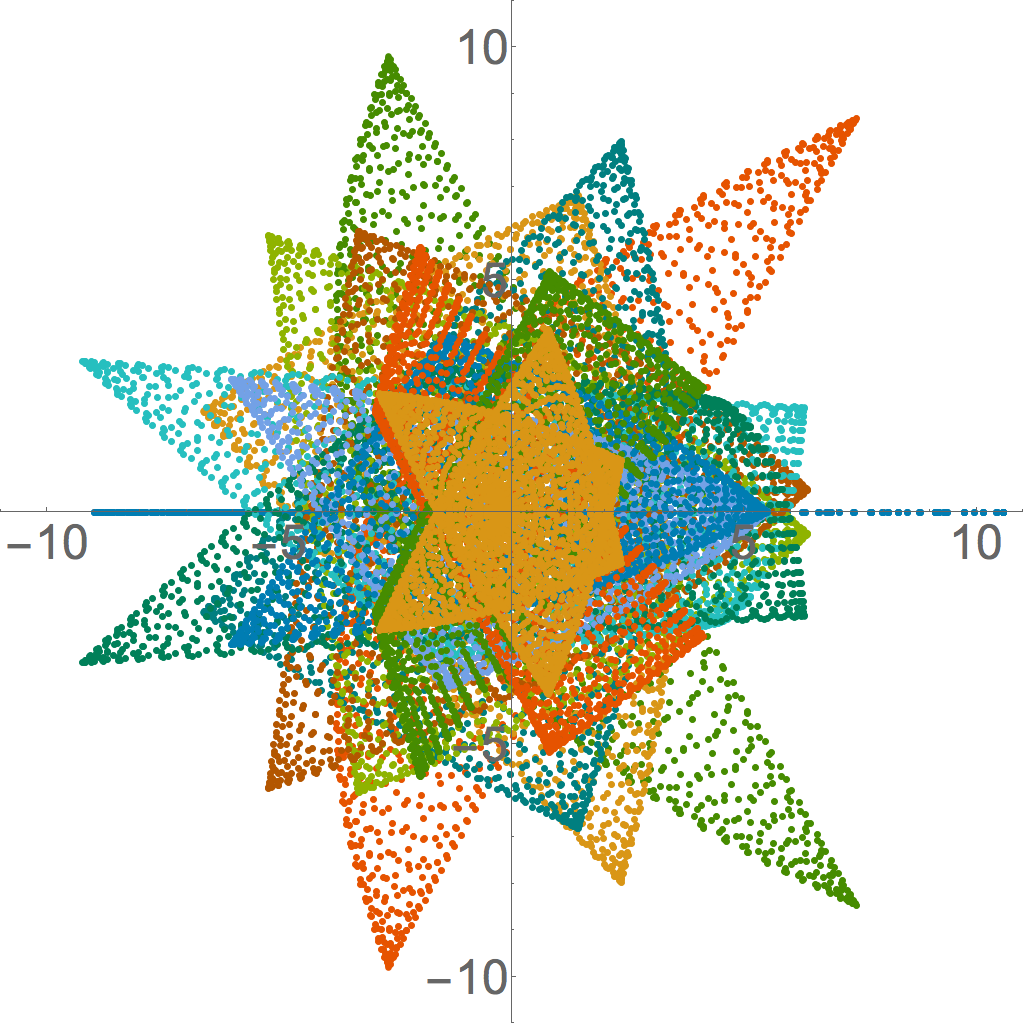}
		\caption{\scriptsize $n=13\cdot 127\cdot 199$, $X=\omega=9247$, $c=127$}
	\end{subfigure}
	\quad
	\begin{subfigure}[b]{0.484\textwidth}
		\includegraphics[width=\textwidth]{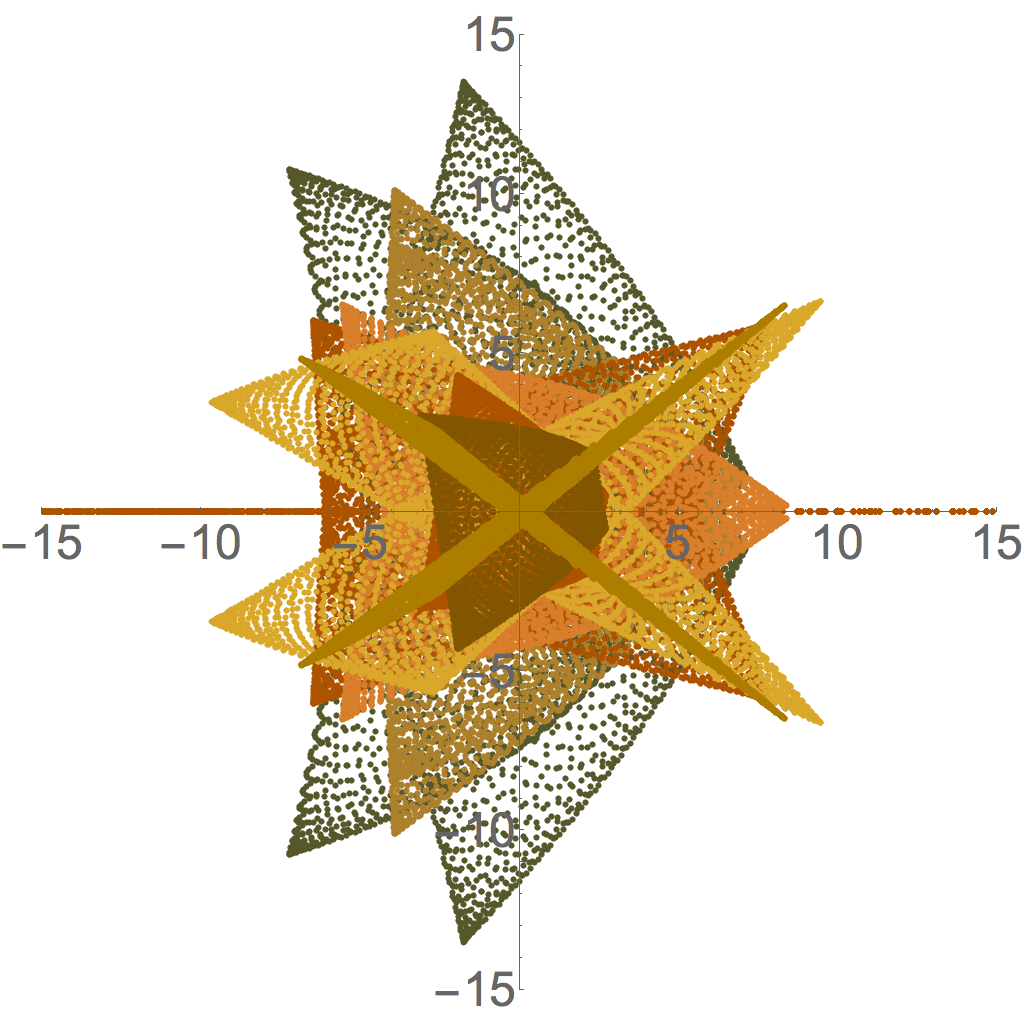}
		\caption{\scriptsize $n=3\cdot 7\cdot 211\cdot 223$, $\omega=710216$, $c=211$}
	\end{subfigure}
	\caption{Mite and moth --- Images of \emph{cyclic supercharacters} correspond to sets of Gaussian periods. For notation and terminology, see Section~\ref{sec:cycsup}.}
	\label{fig:mite}
\end{figure}

The main ingredient in Gauss' argument is an exponential sum known as a \emph{Gaussian period}. Denoting the cardinality of a set $S$ by $|S|$, if $p$ is an odd prime number and $\omega$ has order $d$ in the unit group $(\Z/p\Z)^\times$, then the \emph{$d$-nomial} Gaussian periods modulo $p$ are the complex numbers
\[
	\sum_{j=0}^{d-1} e\(\frac{\omega^j y}{p}\),
\]
where $y$ belongs to $\Z/p\Z$ and $e(\theta)$ denotes $\operatorname{exp}(2\pi i\theta)$ for real $\theta$. Following its appearance in \emph{Disquisitiones}, Gauss' cyclotomy drew the attention of other mathematicians who saw its potential use in their own work. In 1879, J.~J.~Sylvester wrote that ``[c]yclotomy is to be regarded \dots\ as the natural and inherent centre and core of the arithmetic of the future''~\cite{Sylvester}. Two of Kummer's most significant achievements depended critically on his study of Gaussian periods: Gauss' work laid the foundation for the proof of Fermat's Last Theorem in the case of regular primes, and later for Kummer's celebrated reciprocity law.

This success inspired Kummer to generalize Gaussian periods in~\cite{Kummer} to the case of composite moduli. Essential to his work was a study of the polynomial $x^d-1$ by his former student, L.~Kronecker, whom Kummer continued to advocate for the better part of both men's careers~\cite{James}. Just as Gaussian periods for prime moduli had given rise to various families of difference sets~\cite{Baumert}, Kummer's composite cyclotomy has been used to explain certain difference sets arising in finite projective geometry~\cite{Lehmers1}. Shortly after Kummer's publication, L.~Fuchs presented a result in~\cite{Fuchs} concerning the vanishing of Kummer's periods that has appeared in several applications by K.~Mahler~\cite{Mahler1,Mahler2}. A modern treatment of Fuchs' result and a further generalization of Gaussian periods can be found in~\cite{Evans1}.

For a positive integer $n$ and positive divisor $d$ of $\phi(n)$, Kummer ``defined'' a $d$-nomial period modulo $n$ to be the sum
\begin{equation}
	\sum_{j=0}^{d-1} e\(\frac{\omega^j y}{n}\),
	\label{eq:baddef}
\end{equation}
where $\omega$ has order $d$ in the unit group $(\Z/n\Z)^\times$ and $y$ ranges over $\Z/n\Z$.  
Unlike the case of prime moduli, however,
there is no guarantee that a generator $\omega$ of order $d$ will exist, or that a subgroup of order $d$ will be unique. 
For example, consider $(\Z/8\Z)^\times$, which contains no element of order $4$, as well as $3$ distinct subgroups of order $2$.
A similar lack of specificity pervaded some of Kummer's other definitions, including his introduction of ideal prime factors, used to prove a weak form of prime factorization for cyclotomic integers. According to H.~M.~Edwards, instead of revealing deficiencies in Kummer's work, these examples suggest ``the mathematical culture \dots\ as Kummer saw it''~\cite{Edwards}.

Fortunately, the ambiguity in Kummer's definition is easily resolved. For $n$ as above and an element $\omega$ of $(\Z/n\Z)^\times$, we define the Gaussian periods generated by $\omega$ modulo $n$ to be the sums in~\eqref{eq:baddef}, where $d$ is the order of $\omega$ and $y$ ranges over $\Z/n\Z$, as before. These periods are closely related to Gauss sums, another type of exponential sum~\cite{Berndt}.

\begin{figure}[ht]
	\begin{subfigure}[b]{0.484\textwidth}
		\includegraphics[width=\textwidth]{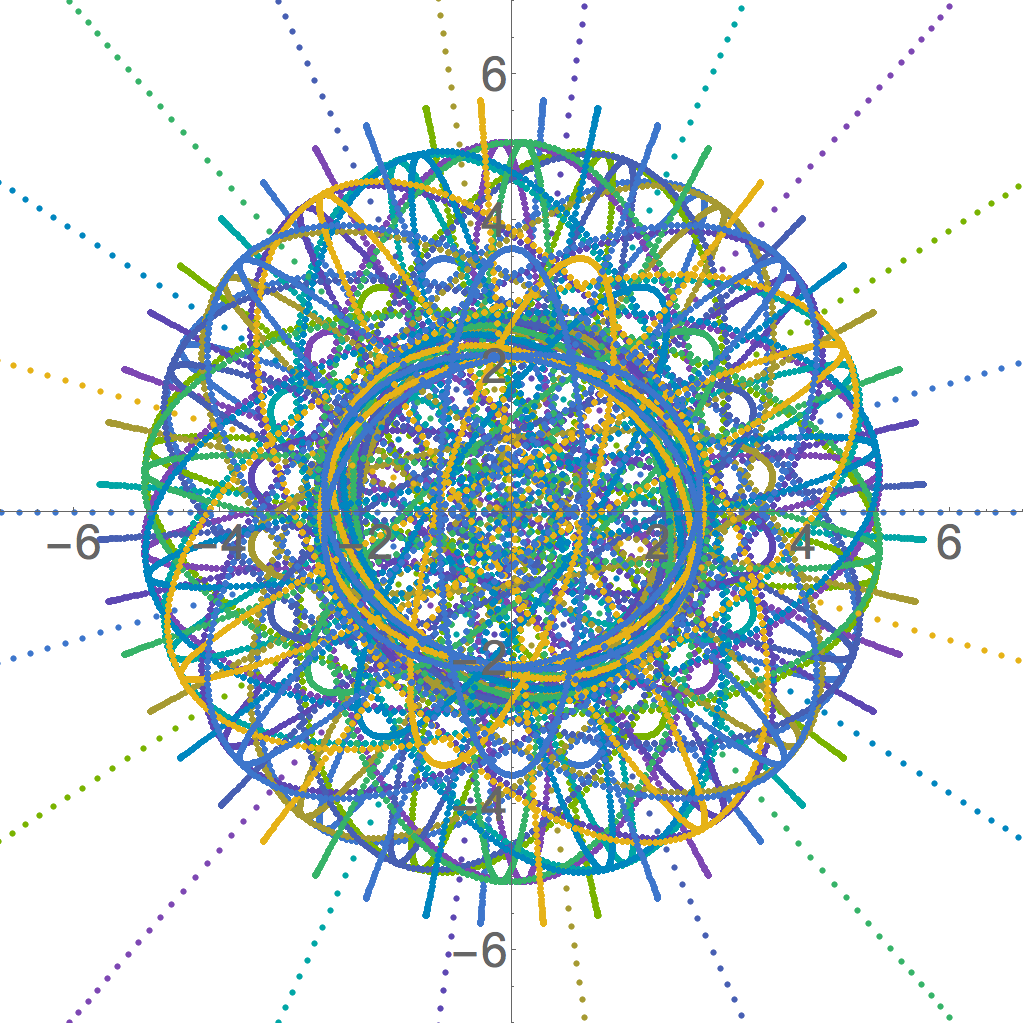}
		\caption{\scriptsize $n=3\cdot 5\cdot 7\cdot 11\cdot 13\cdot 17$, $\omega = 254$, $c=11$}
		\label{fig:atom1}
	\end{subfigure}
	\quad
	\begin{subfigure}[b]{0.484\textwidth}
		\includegraphics[width=\textwidth]{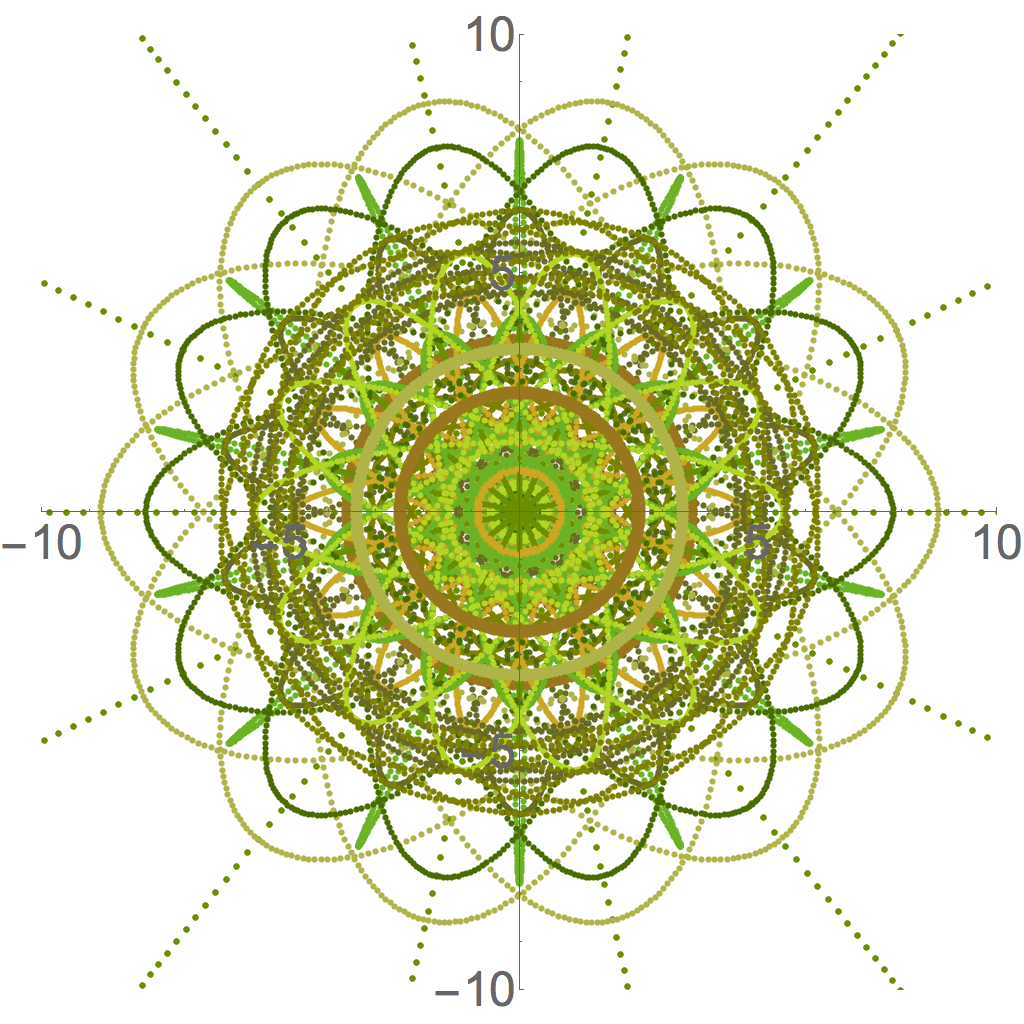}
		\caption{\scriptsize $n=3^2\cdot 5^2\cdot 7 \cdot 17^2$, $\omega=3599$, $c=17^2$}
		\label{fig:atom2}
	\end{subfigure}
	\caption{Atoms --- Images of \emph{cyclic supercharacters} correspond to sets of Gaussian periods. For notation and terminology, see Section~\ref{sec:cycsup}.}
	\label{fig:atom}
\end{figure}


\section{Cyclic supercharacters}
\label{sec:cycsup}

In 2008, P.~Diaconis and I.~M.~Isaacs introduced the theory of supercharacters axiomatically~\cite{DiaconisIsaacs}, building upon seminal work of C.~Andr\'{e} on the representation theory of unipotent matrix groups~\cite{An95, An01}. Supercharacter techniques have been used to study the Hopf algebra of symmetric functions of noncommuting variables~\cite{Aguiar}, random walks on upper triangular matrices~\cite{ACDiSt04}, combinatorial properties of Schur rings~\cite{DiTh09, ThVe09, Th10}, and Ramanujan sums~\cite{Fowler}.

To make an important definition, we divert briefly to the character theory of finite groups. Let $G$ be a finite group with identity 0, $\mathcal{K}$ a partition of $G$, and $\mathcal{X}$ a partition of the set of irreducible characters of $G$. The ordered pair $(\mathcal{X},\mathcal{K})$ is called a \emph{supercharacter theory} for $G$ if $\{0\}\in\mathcal{K}$, $|\mathcal{X}|=|\mathcal{K}|$, and for each $X\in\mathcal{X}$, the function
\[
	\sigma_X=\sum_{\chi\in X} \chi(0)\chi
\]
is constant on each $K\in\mathcal{K}$. The functions $\sigma_X:G\to\C$ are called \emph{supercharacters} and the elements of $\mathcal{K}$ are called \emph{superclasses}.
	
Since $\Z/n\Z$ is abelian, its irreducible characters are group homomorphisms $\Z/n\Z\to \C^\times$. Namely, for each $x$ in $\Z/n\Z$, there is an irreducible character $\chi_x$ of $\Z/n\Z$ given by $\chi_x(y)=e(\frac{xy}{n})$. For a subgroup $\Gamma$ of $(\Z/n\Z)^{\times}$, let $\mathcal{K}$ denote the partition of $\Z/n\Z$ arising from the action $a\cdot x = ax$ of $\Gamma$. The action $a \cdot \chi_x = \chi_{a^{-1}x}$ of $\Gamma$ on the irreducible characters of $\Z/n\Z$ yields a compatible partition $\mathcal{X}$ making $(\mathcal{X},\mathcal{K})$ a supercharacter theory on $\Z/n\Z$. The corresponding supercharacters are 
\begin{equation}
	\sigma_X(y) = \sum_{x \in X} e\left( \frac{xy}{n} \right).
	\label{eq:Sigma}
\end{equation}

For a positive integer $n$ and an orbit $X$ of $\Z/n\Z$ under the multiplication action of a cyclic subgroup $\langle \omega \rangle$ of  $(\Z/n\Z)^\times$, we define the \emph{cyclic supercharacter} $\sigma_X:\Z/n\Z\to \C$ by~\eqref{eq:Sigma}.  The values of these functions 
are Gaussian periods in the sense of Kummer~\cite{DukeGarciaLutz}. For applications of supercharacter theory to other exponential sums, see~\cite{Brumbaugh1, Brumbaugh2, Fowler}.

\begin{figure}[b]
	\begin{subfigure}[b]{0.484\textwidth}
		\includegraphics[width=\textwidth]{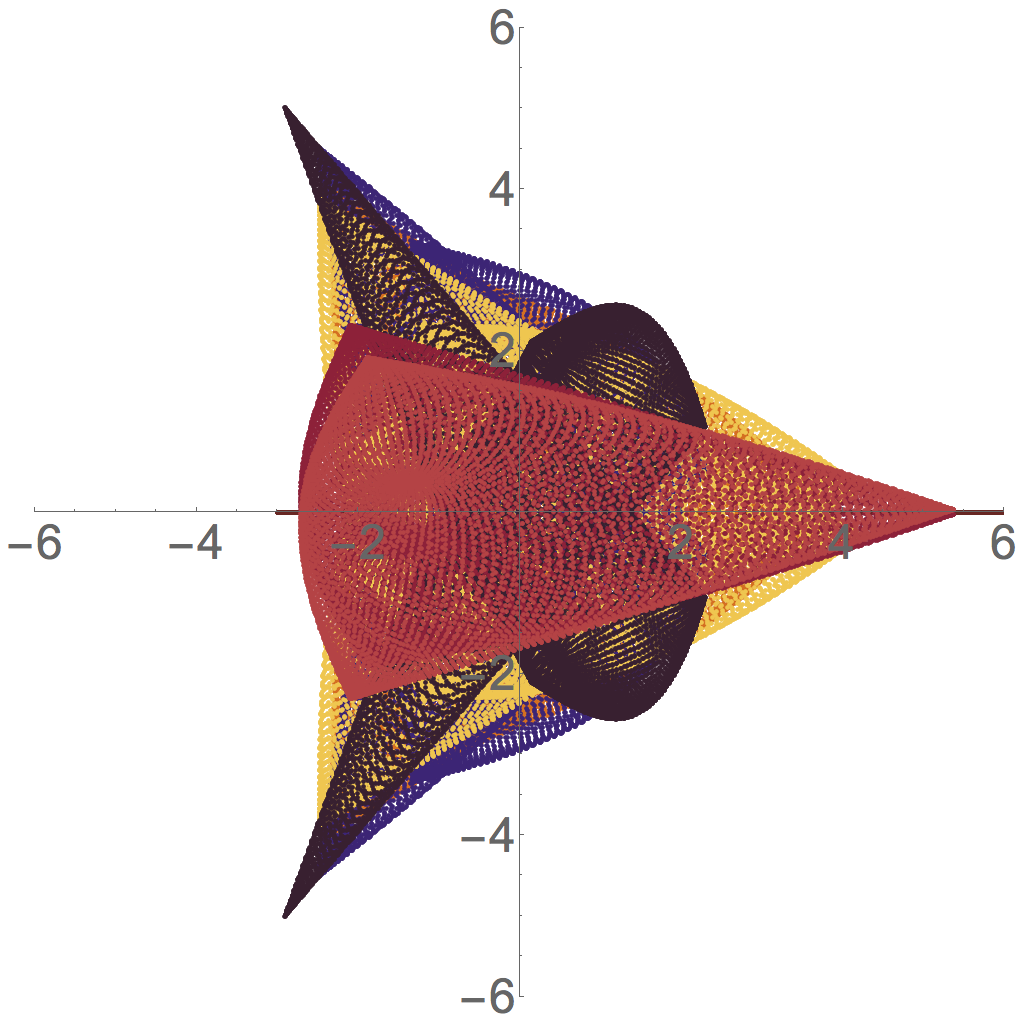}
		\caption{\scriptsize $n=31\cdot 73\cdot 211$, $\omega=2547$, $c=31$}
	\end{subfigure}
	\quad
	\begin{subfigure}[b]{0.484\textwidth}
		\includegraphics[width=\textwidth]{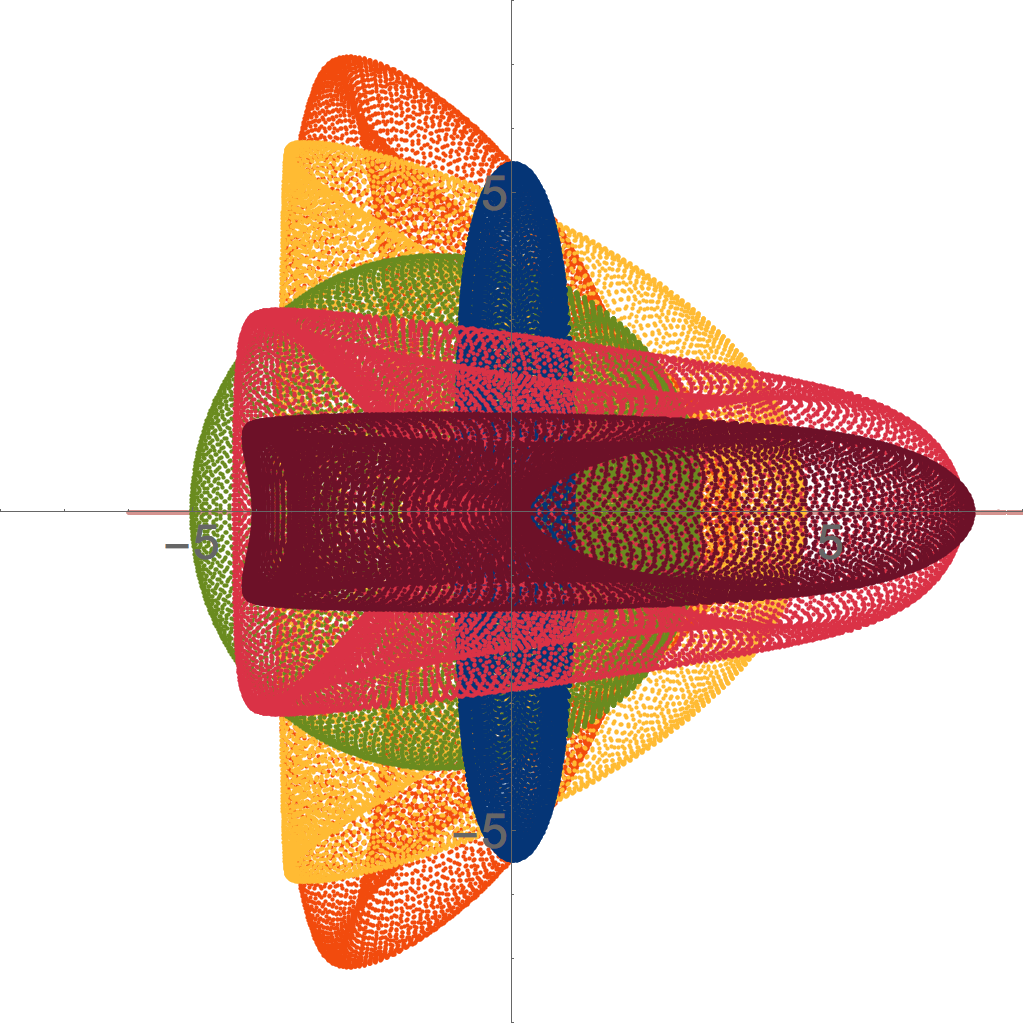}
		\caption{\scriptsize $n=3\cdot 31\cdot 73\cdot 211$, $\omega= 1463$, $c=73$}
	\end{subfigure}
	\caption{Bird and spacecraft --- images of cyclic supercharacters.}
	\label{fig:ship}
\end{figure}

We are now in a position to clarify the captions and colorizations of the numerous figures. Unless specified otherwise, the image appearing in each figure is the image in $\C$ of the cyclic supercharacter $\sigma_{\inner{\omega}} :\Z/n\Z\to \C$, where $\omega$ belongs to $(\Z/n\Z)^\times$, and $\inner{\omega}=\inner{\omega}1$ denotes the orbit of 1 under the action of the subgroup generated by $\omega$. Conveniently, the image of \emph{any} cyclic supercharacter is a scaled subset of the image of one having the form $\sigma_{\inner{\omega}}$~\cite[Proposition 2.2]{DukeGarciaLutz}, so a restriction of our attention to orbits of 1 is natural. Moreover, the image of $\sigma_{\inner{\omega}}$ is the set of Gaussian periods generated by $\omega$ modulo $n$, bringing classical relevance to these figures.

To colorize each image, we fix a proper divisor $c$ of $n$ and assign a color to each of the \emph{layers}
\[
	\big\{\sigma_{\langle \omega\rangle 1}(y)\mid y\equiv j\!\!\!\!\pmod{c}\big\},
\]
for $j=0,1,\ldots,c-1$.  Different choices of $c$ result in different ``layerings.'' For many images, certain values of $c$ yield colorizations that separate distinct graphical components.

Predictable layering occurs when the image of a cyclic supercharacter contains several rotated copies of a proper subset. We say that a subset of $\C$ has \emph{$k$-fold dihedral symmetry} if it is invariant under complex conjugation and rotation by $\frac{2\pi}{k}$ about the origin. For example, the image pictured in Figure~\ref{fig:atom1} has 11-fold dihedral symmetry, while the symmetry in Figure~\ref{fig:atom2} is 7-fold. The image of a cyclic supercharacter $\sigma_{\inner{\omega}}:\Z/n\Z\to \C$ has $k$-fold dihedral symmetry, where $k=\gcd\(n,\omega-1\)$~\cite[Proposition 3.1]{DukeGarciaLutz}. In this situation, taking $c=k$ results in exactly $k$ layers that are rotated copies of one another.

\begin{figure}[ht]
	\begin{subfigure}[b]{0.484\textwidth}
		\includegraphics[width=\textwidth]{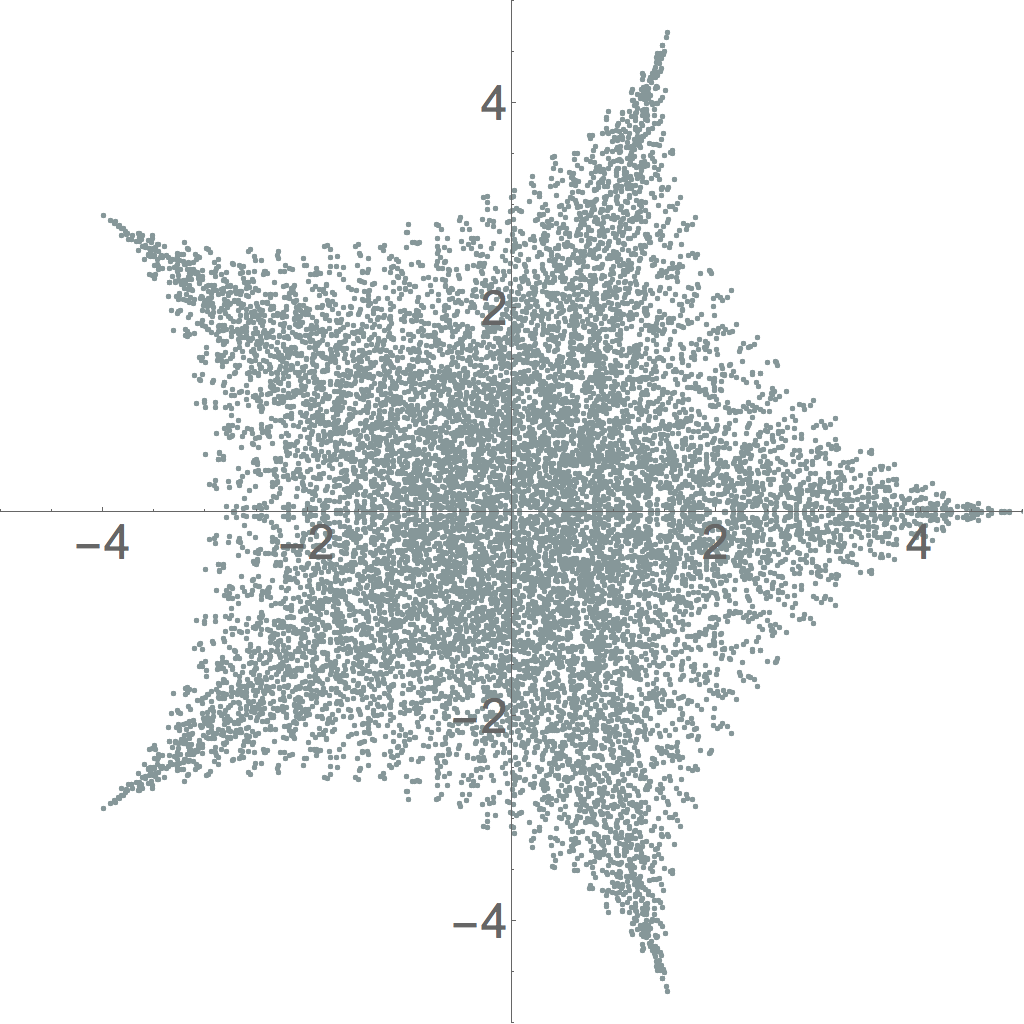}
		\caption{\scriptsize $m=251\cdot281$, $\omega_m=54184$}
	\end{subfigure}
	\quad
	\begin{subfigure}[b]{0.484\textwidth}
		\includegraphics[width=\textwidth]{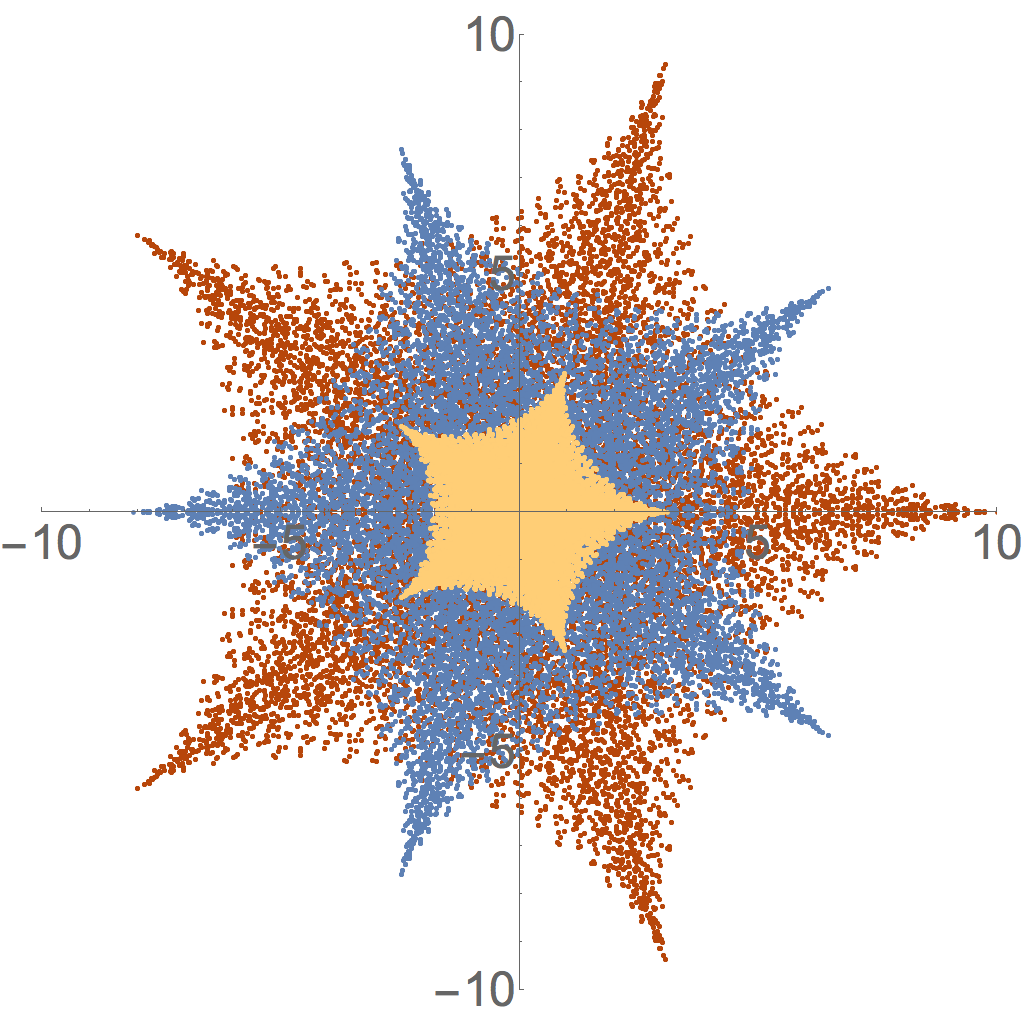}
		\caption{\scriptsize $mn=5\cdot 251\cdot 281$, $\omega=54184$, $c=5$}
	\end{subfigure}
	\caption{The image on the right is the product set of the image on the left and the image $\{2,\frac{1}{2}(\pm\sqrt{5}-1)\}$,
	as in~\eqref{eq:mult}.}
	\label{fig:product}
\end{figure}

\begin{figure}[ht]
	\includegraphics[width=0.7\textwidth]{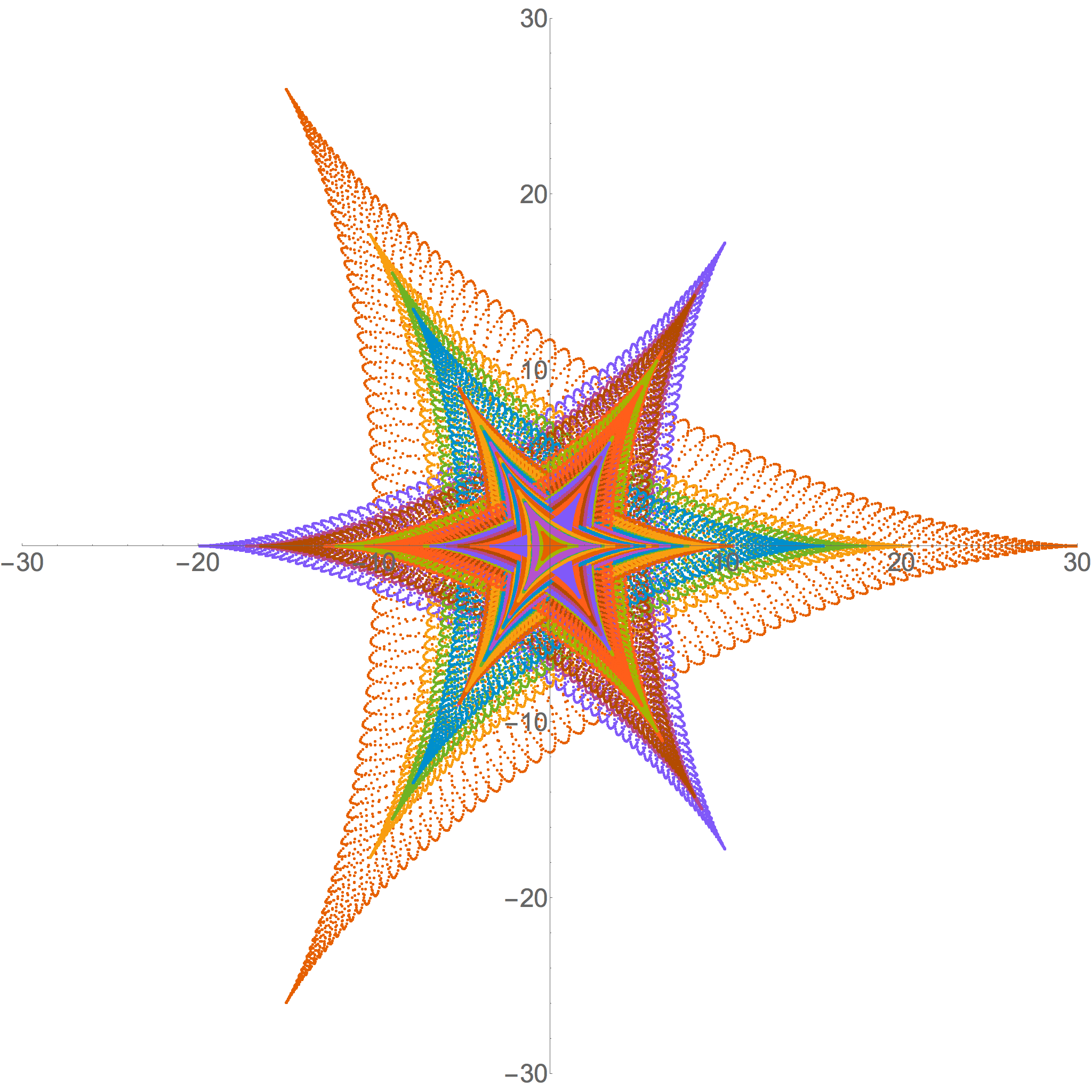}
	\caption{$n=127^2\cdot 401$, $\omega=6085605$, $c=401$}
	\label{fig:nested}
\end{figure}

In addition to the behaviors above, certain cyclic supercharacters enjoy a multiplicative property \cite[Theorem 2.1]{DukeGarciaLutz}. Specifically, if $\gcd(m,n)=1$ and $\omega\mapsto (\omega_m,\omega_n)$ under the isomorphism $(\Z/mn\Z)^\times\to (\Z/m\Z)^\times \times(\Z/n\Z)^\times$ afforded by the Chinese Remainder Theorem, where the multiplicative orders of $\omega_m$ and $\omega_n$ are coprime, then
\begin{equation}
	\sigma_{\langle \omega\rangle }(\Z/mn\Z)=
	\big\{wz\in\C:(w,z)\in \sigma_{\langle \omega_m\rangle }(\Z/m\Z)
	\times \sigma_{\langle \omega_n\rangle }(\Z/n\Z)\big\}.
\label{eq:mult}
\end{equation}
This can be used to explain the images of cyclic supercharacters featuring ``nested'' copies of a given shape. For examples of this phenomenon, see Figures~\ref{fig:product} and~\ref{fig:nested}.

\section{Asymptotic behavior}
\label{sec:asym}
	
	In this section, we restrict our attention to cyclic supercharacters $\sigma_X:\Z/q\Z\to\C$, where $q$ is a power of an odd prime and $X=\inner{\omega}$ is an orbit of 1. The Gaussian periods attained as values of these supercharacters have been applied in various settings \cite{BaumertMillsWard,Bayer,Hoshi}. Plotting the functions $\sigma_X$ in this case reveals asymptotic patterns that have, until recently, gone unseen. Before proceeding, we recall several definitions and results.

\subsection{Uniform distribution mod 1}

Let $m$ be a positive integer and $\Lambda$ a finite subset of $\R^m$. We write
\[
	\widehat{\Lambda}=\big\{(\lambda_1-\floor{\lambda_1},\ldots,\lambda_m-\floor{\lambda_m}) \in [0,1)^m : (\lambda_1,\ldots,\lambda_m)\in \Lambda\big\},
\]
where $\floor{\cdot}$ denotes the greatest integer function. The \emph{discrepancy} of the set $\Lambda$ is
\[
	\sup_B \left| \frac{|B\cap \widehat{\Lambda}|}{|\widehat{\Lambda}|}-\operatorname{vol}(B)\right|,
\]
where the supremum is taken over all boxes $B=[a_1,b_1)\times \cdots\times [a_m,b_m)\subset [0,1)^m$ and $\operatorname{vol}(B)$ denotes the volume of $B$. We
say that a sequence $(\Lambda_n)_{n=1}^\infty$ of finite subsets of $\R^m$ is \emph{uniformly distributed mod 1} if the discrepancy of $\Lambda_n$ 
tends to zero as $n\to\infty$.

Many accessible examples exist in the case $m=1$. For instance, a well-known theorem of Weyl states that if $\alpha$ is an irrational number, then the sequence
\begin{equation}
	\big(\{k\alpha\in \R : k=0,1,\ldots,n-1\} \big)_{n=1}^\infty
	\label{eq:equidist}
\end{equation}
is uniformly distributed mod 1 \cite{Weyl2}. A lovely result of Vinogradov is that the same holds when $k$ above is replaced by the $k$th prime number \cite{Vinogradov}. While it is known that the sequence
\[
	\big(\{\theta,\theta^2,\ldots,\theta^n\} \big)_{n=1}^\infty
\]
is uniformly distributed mod 1 for almost every $\theta>1$, specific cases such as $\theta=\frac{3}{2}$ are not well understood \cite{Elkies}. Surprisingly perhaps, the sequence
\[
	\big(\{\log 1,\log 2,\ldots,\log n\}\big)_{n=1}^\infty
\]
is not uniformly distributed mod 1. This example and several others are elaborated in \cite{Kuipers} using the following crucial characterization, also due to Weyl \cite{Weyl}.

\begin{Lemma}[Weyl's criterion]\label{Lemma:Weyl}
A sequence $(\Lambda_n)_{n=1}^\infty$ of finite subsets of $\R^m$ is uniformly distributed mod 1 if and only if for each $\vec{v}$ in $\Z^m$ we have
\[
	\lim_{n\to \infty} \frac{1}{|\Lambda_n|}
	\sum_{\vec{u}\in \Lambda_n} e(\vec{u}\cdot\vec{v})=0.
\]
\end{Lemma}

For example, let $(\Lambda_n)_{n=1}^\infty$ be the sequence in \eqref{eq:equidist}. Since $\alpha$ is irrational,
\[
	\frac{1}{|\Lambda_n|}\sum_{u\in \Lambda_n} e(uv)
	= \frac{1}{n}\sum_{k=0}^{n-1} e(v\alpha)^k
	= \frac{1}{n}\left(\frac{1-e(v\alpha)^n}{1-e(v\alpha)}\right)
\]
for each nonzero $v\in\Z$.  Consequently,
\[
	\left|\frac{1}{n}\left(\frac{1-e(v\alpha)^n}{1-e(v\alpha)}\right)\right|
	\leq \frac{1}{n}\frac{2}{|1-e(v\alpha)|}\to 0
\]
as $n\to \infty$, so Lemma \ref{Lemma:Weyl} confirms that \eqref{eq:equidist} is uniformly distributed mod 1.

\subsection{Cyclotomic polynomials}

For a positive integer $d$, the \emph{$d$th cyclotomic polynomial} $\Phi_d(x)$ is defined by the formula
\[
	\Phi_d(x)=
	\prod_{\substack{k=1,2,\ldots,d\\\gcd(k,d)=1}}
	\(x-e\(\tfrac{k}{d}\)\).
\]
It can be shown that $\Phi_d(x)$ is of degree $\phi(d)$ and belongs to $\Z[x]$. In \emph{Disquisitiones}, Gauss showed that $\Phi_d(x)$ is irreducible, hence the minimal polynomial of any primitive $d$th root of unity, over $\Q$. The first several cyclotomic polynomials are
\begin{align*}
\Phi_1(t)&=x-1,\\
\Phi_2(t)&=x+1,\\
\Phi_3(t)&=x^2+x+1,\\
\Phi_4(t)&=x^2+1,\\
\Phi_5(t)&=x^4+x^3+x^2+x+1.
\end{align*}
In these examples, the coefficients have absolute value at most $1$. In 1938, N.~G.~Chebotar\"{e}v asked  whether this phenomenon continues for all factors of $x^d-1$ and all values of $d$ \cite{Grassl}. Three years later, V.~Ivanov showed that while the pattern holds for $d<105$, one coefficient of $\Phi_{105}(x)$ is $-2$.  Unbeknownst to either mathematician,
A.~S.~Bang had solved Chebotar\"{e}v's challenge more than forty years earlier \cite{Bloom}.

\subsection{Main Theorem}

We require a lemma essentially due to G.~Myerson. The original aim of the result was to count the number of ways to write an arbitrary element of $(\Z/q\Z)^\times$ as a sum of elements, one in each coset of a fixed subgroup of $(\Z/q\Z)^\times$. For our purposes, Myerson's lemma will be critical in discussing the asymptotic behavior of the images of cyclic supercharacters. Throughout, we let $\omega_q$ denote a primitive $d$th root of unity in $\Z/q\Z$, in which $q=p^a$ is a power of an odd prime $p$, and
\begin{equation}
	\Lambda_q=
	\left\{\frac{\ell}{q}\(1,\omega_q,\omega_q^2,\ldots,\omega_q^{\phi(d)-1}\)
	\in [0,1)^{\phi(d)} : \ell =0,1,\ldots,q-1\right\}.
	\label{eq:omegaset}
\end{equation}

\begin{figure}[ht]
	\begin{subfigure}[b]{0.3\textwidth}
		\includegraphics[width=\textwidth]{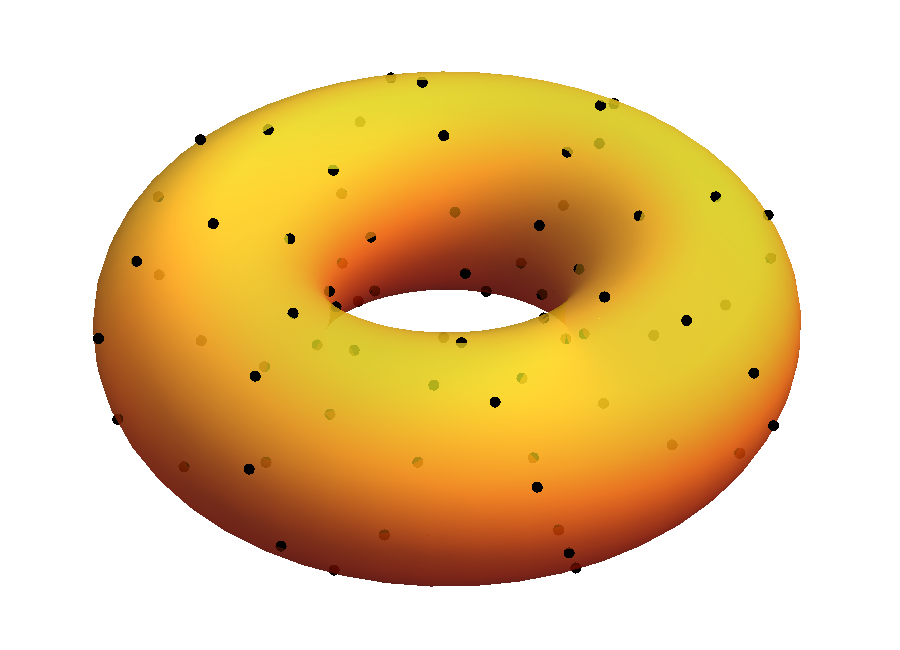}
		\caption{\scriptsize $q=73$}
	\end{subfigure}
	\quad
	\begin{subfigure}[b]{0.3\textwidth}
		\includegraphics[width=\textwidth]{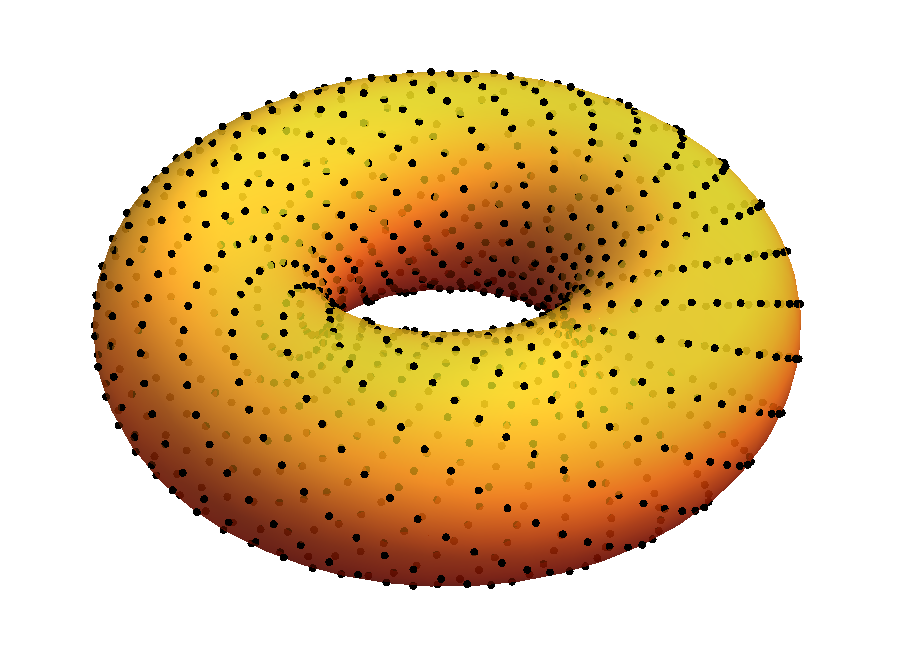}
		\caption{\scriptsize $q=31^2$}
	\end{subfigure}
	\quad
	\begin{subfigure}[b]{0.3\textwidth}
		\includegraphics[width=\textwidth]{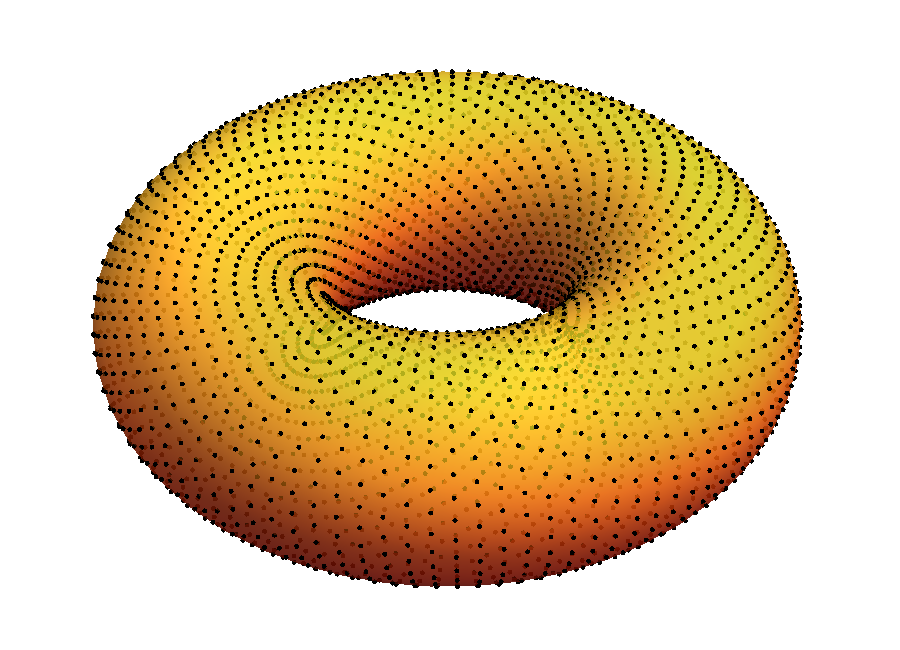}
		\caption{\scriptsize $q=3571$}
	\end{subfigure}
	\caption{The sequence $(\Lambda_q)_{q\equiv 1\!\!\pmod{3}}$, defined in \eqref{eq:omegaset}, is uniformly distributed mod 1. Here we plot several of the sets $\Lambda_q\subset [0,1)^2$, identifying $[0,1)^2$ with a torus.}
\end{figure}

\begin{Lemma}[Myerson \cite{Myerson}]
The sequence $(\Lambda_q)_{q\equiv 1\!\!\pmod{d}}$ is uniformly distributed mod 1.
\begin{proof}
Let $\vec{v}=(v_0,\ldots,v_{\phi(d)-1})$ be nonzero in $\Z^{\phi(d)}$ and let $f\in\Z[x]$ be given by
\[
f(x)=v_0+v_1x+\cdots+v_{\phi(d)-1}x^{\phi(d)-1}.
\]
Writing $r=\frac{q}{\gcd(q,f(\omega_q))}$, we notice that
\begin{align}
	\sum_{\vec{u}\in \Lambda_q} e(\vec{u}\cdot\vec{v})
		&=\sum_{\ell=0}^{q-1} e\(\frac{\ell f(\omega_q)}{q}\)\nonumber\\
		&=\sum_{k=0}^{q/r-1} \sum_{j=0}^{r-1} e\(\frac{(kr+j)f(\omega_q)}{q}\)\nonumber\\
		&=\sum_{k=0}^{q/r-1} \sum_{j=0}^{r-1} e\(\frac{kf(\omega_q)}{\gcd(q,f(\omega_q))}+\frac{jf(\omega_q)}{q}\)\nonumber\\
		&=\sum_{k=0}^{q/r-1} \sum_{j=0}^{r-1} e\(\frac{jf(\omega_q)}{q}\)\nonumber\\
		&=\frac{q}{r}\sum_{j=0}^{r-1} e\(\frac{jf(\omega_q)}{q}\)\nonumber\\
		&=\begin{cases}q&\mbox{if } q|f(\omega_q),\\ 0&\mbox{otherwise}.\\\end{cases}
	\label{eq:qcomp}
\end{align}
Since $\Phi_d(x)$ is irreducible over $\Q$ and of greater degree than $f(x)$, we see that $\gcd(f(x),\Phi_d(x))=1$ in $\Q[x]$. From this we obtain $a(x)$ and $b(x)$ in $\Z[x]$ such that $a(x)f(x)+b(x)\Phi_d(x)=s$
for some $s\in\Z$. Evaluating at $x=\omega_q$ reveals that $a(\omega_q)f(\omega_q)\equiv s\pmod{q}$, so $q|s$ whenever $q|f(\omega_q)$. Hence $q|f(\omega_q)$ for at most finitely many odd prime powers 
$q\equiv 1\!\!\pmod{d}$. It follows from~\eqref{eq:qcomp} that
\[
	\lim_{\substack{q\to\infty\\ q\equiv1\!\!\!\!\pmod{d}}}
	\frac{1}{|\Lambda_q|} \sum_{\vec{u}\in\Lambda_q} e(\vec{u}\cdot\vec{v})
	= 0.
\]
Appealing to Lemma \ref{Lemma:Weyl}
completes the proof.
\end{proof}
\label{lem:myerson}
\end{Lemma}

The following theorem summarizes our current understanding of the asymptotic behavior of cyclic supercharacters. A special case supplies a geometric description of the set of $|X|$-nomial periods modulo $p$ considered by Gauss, shedding new light on these classical objects. We let $\T$ denote the unit circle in $\C$.

\begin{Theorem}[Duke--Garcia--Lutz]\label{DukeThm}
If $\sigma_X:\Z/q\Z\to\C$ is a cyclic supercharacter, where $q=p^a$ is a power of an odd prime, $X$ is an orbit of 1, and $|X|=d$ divides $p-1$, then the image of $\sigma_X$ is contained in the image of the Laurent polynomial function $g_d:\T^{\phi(d)}\to \C$ defined by
\begin{equation}\label{eq:gfunc}
	g_d(z_1,z_2,\ldots,z_{\phi(d)})=
	\sum_{k=0}^{d-1}\prod_{j=0}^{\phi(d)-1} z_{j+1}^{c_{jk}},
\end{equation}
where the $c_{jk}$ are given by the relation
\begin{equation}
	x^k\equiv \sum_{j=0}^{\phi(d)-1} c_{jk}x^j
	\!\!\!\!\pmod{\Phi_d(x)}.
	\label{eq:cjk}
\end{equation}
Moreover, for a fixed $|X|=d$, as $q\equiv 1\pmod d$ tends to infinity, every nonempty open disk in the image of $g$ eventually contains points in the image $\sigma_X(\Z/q\Z)$.

\begin{proof}
Since the elements $1,\omega_q,\ldots,\omega_q^{\phi(d)-1}$ form a $\Z$-basis for $\Z[\omega_q]$~\cite[p.~60]{Neukirch}, for $k=0,1,\ldots,d-1$ we can write
\[
	\omega_q^k\equiv \sum_{j=0}^{\phi(d)-1} c_{jk} \omega_p^{j} \!\!\!\!\pmod{q},
\]
where the integers $c_{jk}$ are given by~\eqref{eq:cjk}. Letting $X=\inner{\omega_q}$, we see that
\[
	\sigma_X(y)=
	\sum_{x \in X} e\(\frac{xy}{q}\)=
	\sum_{k=0}^{d-1} e\(\frac{\omega_q^ky}{q}\)=
	\sum_{k=0}^{d-1} e\(\sum_{j=0}^{\phi(d)-1} c_{jk}\frac{\omega_q^{j}y}{q}\)=
	\sum_{k=0}^{d-1}\prod_{j=0}^{\phi(d)-1} e\(\frac{\omega_q^{j}y}{q}\)^{c_{jk}},
\]
whence the image of $\sigma_X$ is contained in the image of the function
$g_d:\T^{\phi(d)}\to\C$ defined in~\eqref{eq:gfunc}. The claim about open disks follows immediately from Lemma~\ref{lem:myerson}.
\end{proof}
\end{Theorem}

\begin{figure}[ht]
	\begin{subfigure}[b]{0.3\textwidth}
		\includegraphics[width=\textwidth]{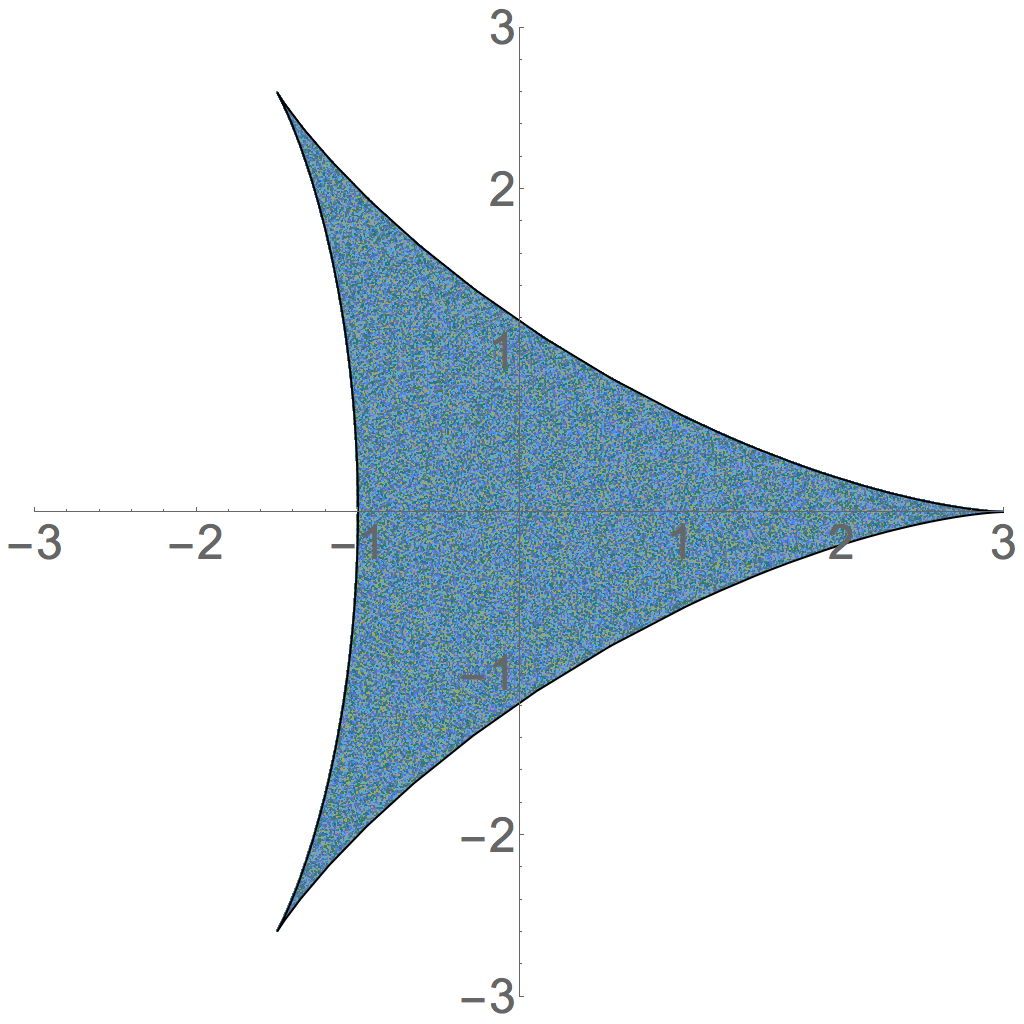}
		\caption{\scriptsize $n=97^3$, $\omega=61074$, $d=3$}
	\end{subfigure}
	\quad
	\begin{subfigure}[b]{0.3\textwidth}
		\includegraphics[width=\textwidth]{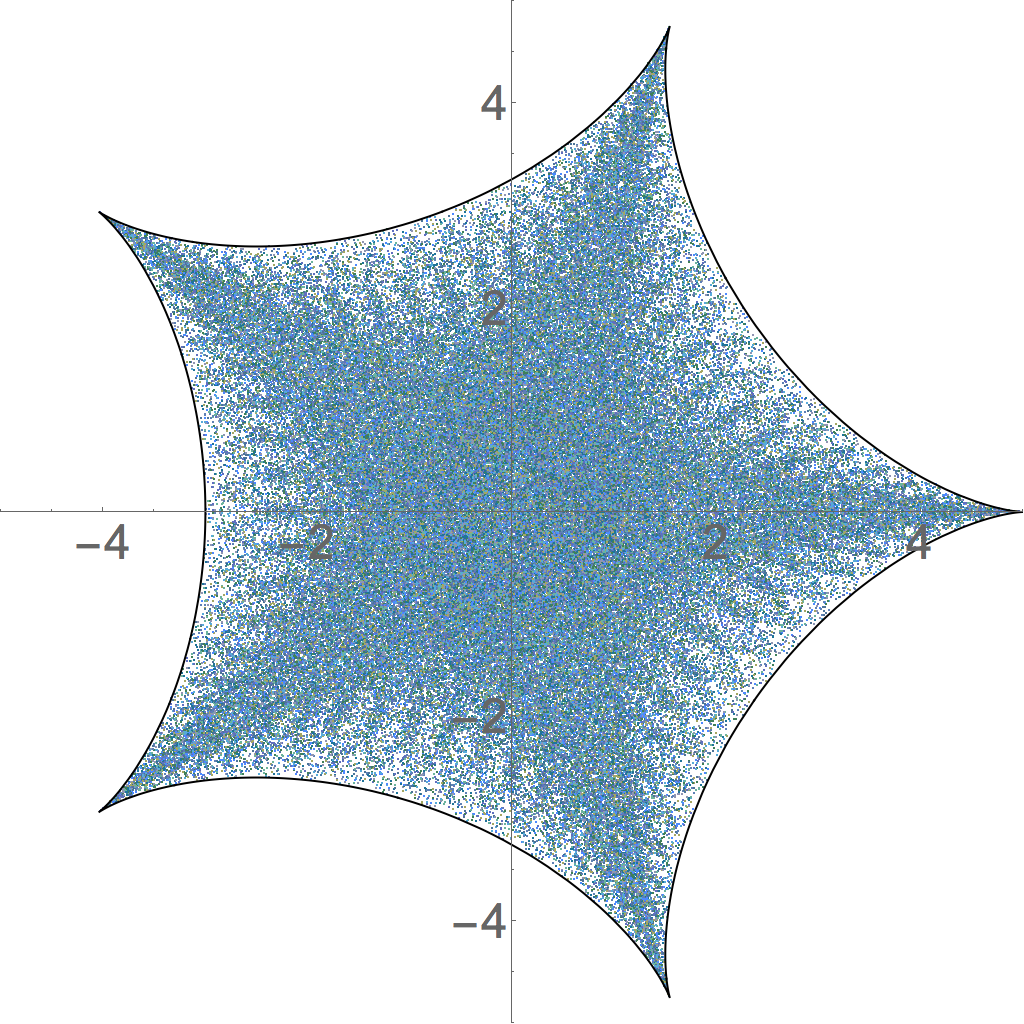}
		\caption{\scriptsize $n=31^4$, $\omega=62996$, $d=5$}
	\end{subfigure}
	\quad
	\begin{subfigure}[b]{0.3\textwidth}
		\includegraphics[width=\textwidth]{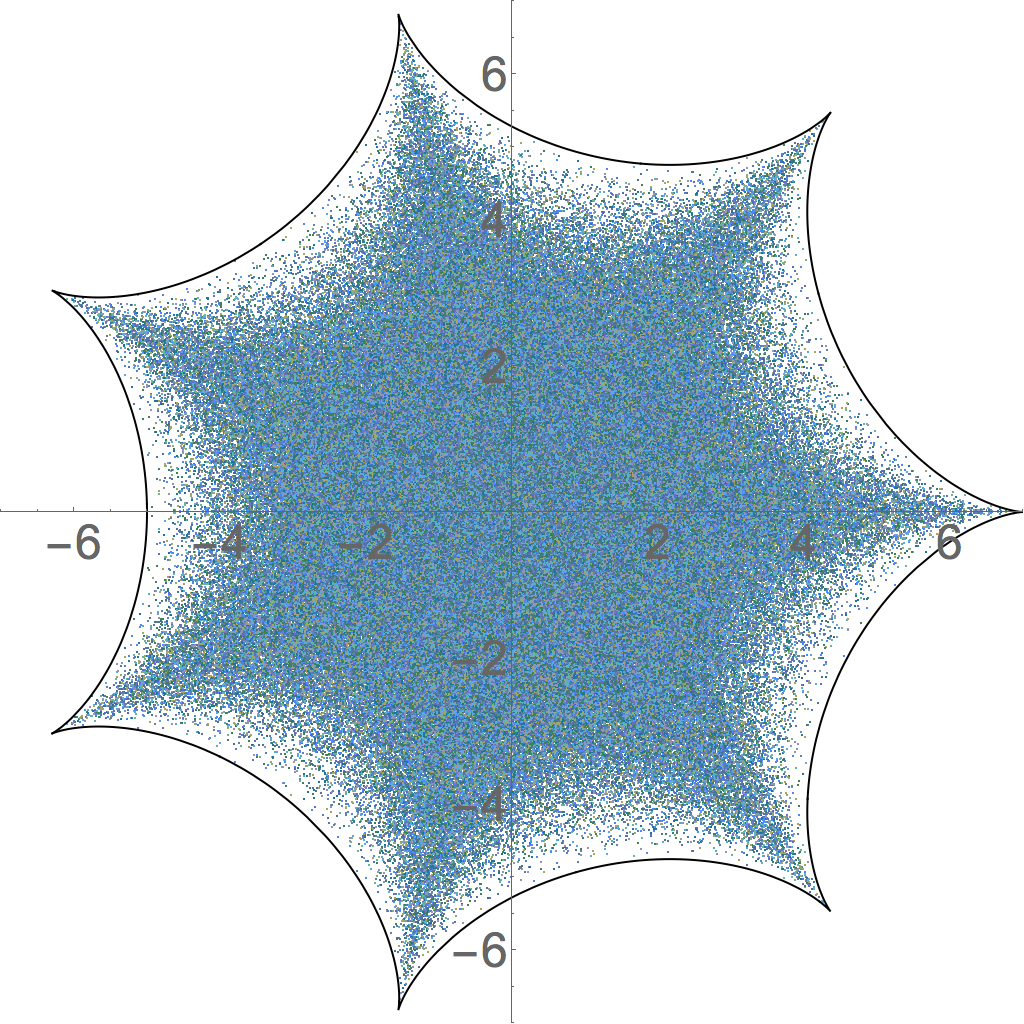}
		\caption{\scriptsize $n=1933^2$, $\omega=537832$, $d=7$}
	\end{subfigure}
	\caption{Certain cyclic supercharacters fill out regions bounded by \emph{hypocycloids}, outlined in black (see Proposition~\ref{prop:hypo}).}
	\label{fig:hypo}
\end{figure}
	
	Several remarks are in order. First, when the hypotheses of Theorem~\ref{DukeThm} are satisfied, we say that $\sigma_X$ \emph{fills out} the image of $g_d$, as illustrated by Figure~\ref{fig:hypo}. The corresponding values of $d$ are given in captions. Second, since every divisor of $\phi(q)=p^{a-1}(p-1)$ is the cardinality of some orbit $X$ under the action of a cyclic subgroup of $(\Z/p^a\Z)^\times$, the requirement that $|X|$ divide $p-1$ might seem restrictive. However, it turns out that if $p$ divides $|X|$, then the image of $\sigma_X$ is equal to a scaled copy of the image of a supercharacter that satisfies the hypotheses of the theorem, except for a single point at the origin \cite[Proposition 2.4]{DukeGarciaLutz}.

\section{Examples}
\label{sec:ex}

As a consequence of Theorem~\ref{DukeThm}, the functions $g_d$ are instrumental in understanding the asymptotic behavior of cyclic supercharacters $\sigma_X:\Z/q\Z\to \C$, where $q$ is a power of an odd prime. Fortunately, whenever the coefficients of $\Phi_d(x)$ are relatively accessible, we can obtain a convenient formula for $g_d$. For example, it is not difficult to show that
\[
	\Phi_{2^b}(x)=x^{2^{b-1}}+1,
\]
for any positive integer $b$. With this, we can compute the integers $c_{jk}$ in~\eqref{eq:cjk} to see that
\[
	g_{2^b}(z_1,z_2,\ldots,z_{2^{b-1}})
	=
	2\sum_{j=1}^{2^{b-1}} \Re (z_j),
\]
where $\Re(z)$ denotes the real part of $z$. Hence the image of $g_{2^b}$ is the real interval $[-2^b,2^b]$. Alternatively, if $r$ is an odd prime, then
\[
	\Phi_{2r}(x)=\sum_{j=0}^{r-1} (-x)^j,
\]
giving
\[
	g_{2r}(z_1,z_2,\ldots,z_{r-1}) = 2\Re\(\frac{z_2z_4\cdots z_{r-1}}{z_1z_3\cdots z_{r-2}}\)+2\sum_{j=1}^{r-1}\Re(z_j).
\]
More generally, $g_d$ is real valued whenever $d$ is even.

\begin{figure}[ht]
	\begin{subfigure}[b]{0.3\textwidth}
		\includegraphics[width=\textwidth]{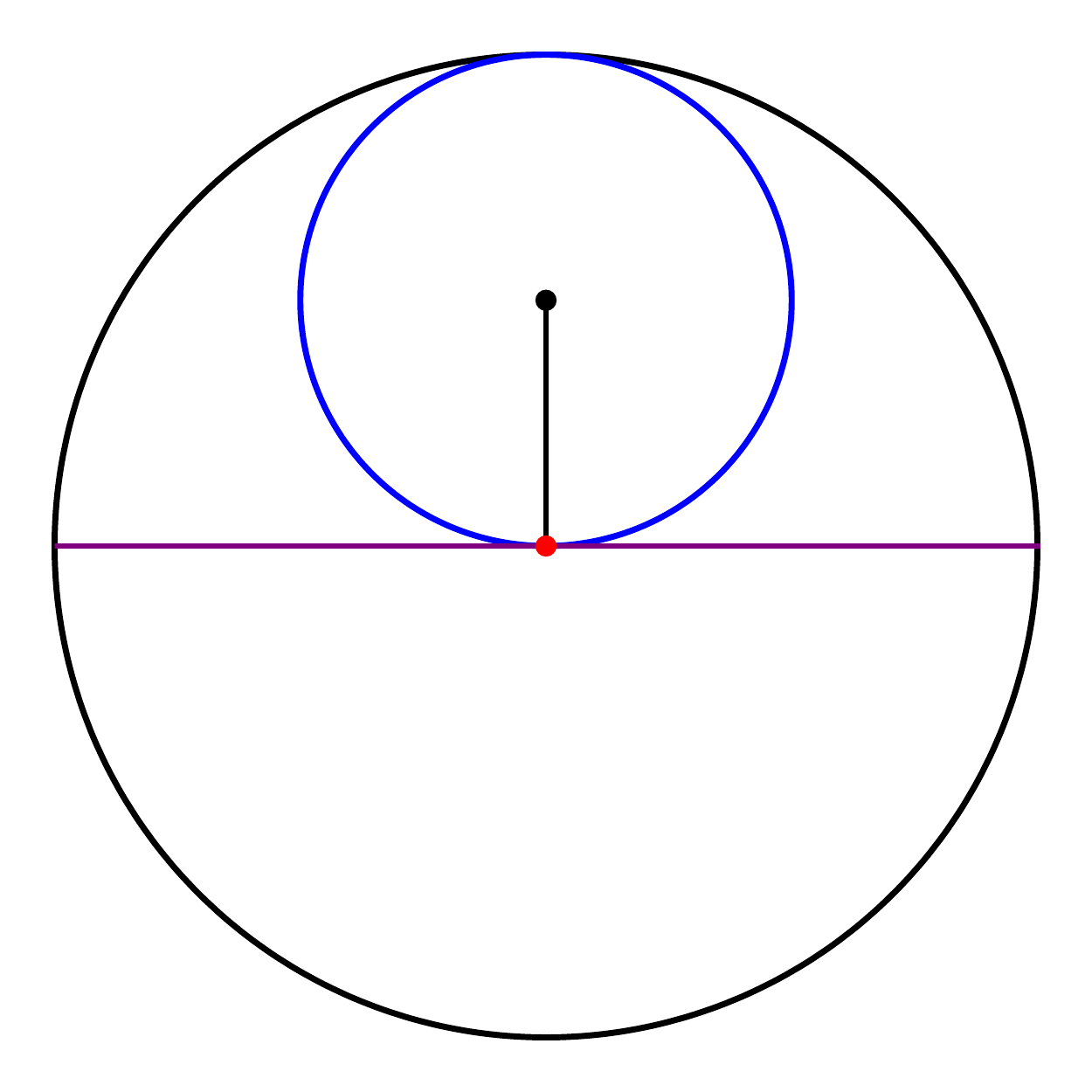}
		\caption{Tusi couple}
	\end{subfigure}
	\quad
	\begin{subfigure}[b]{0.3\textwidth}
		\includegraphics[width=\textwidth]{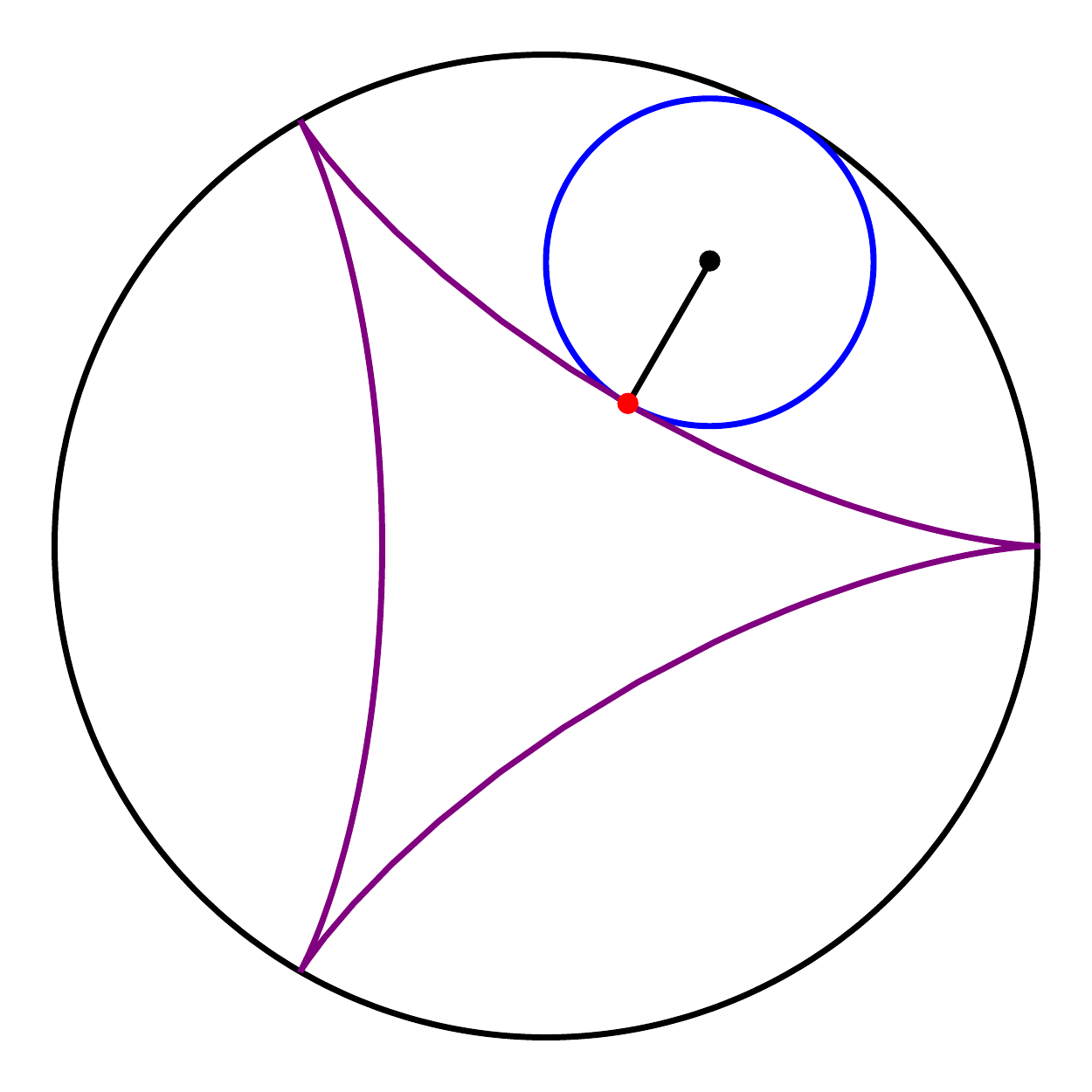}
		\caption{Deltoid}
	\end{subfigure}
	\quad
	\begin{subfigure}[b]{0.3\textwidth}
		\includegraphics[width=\textwidth]{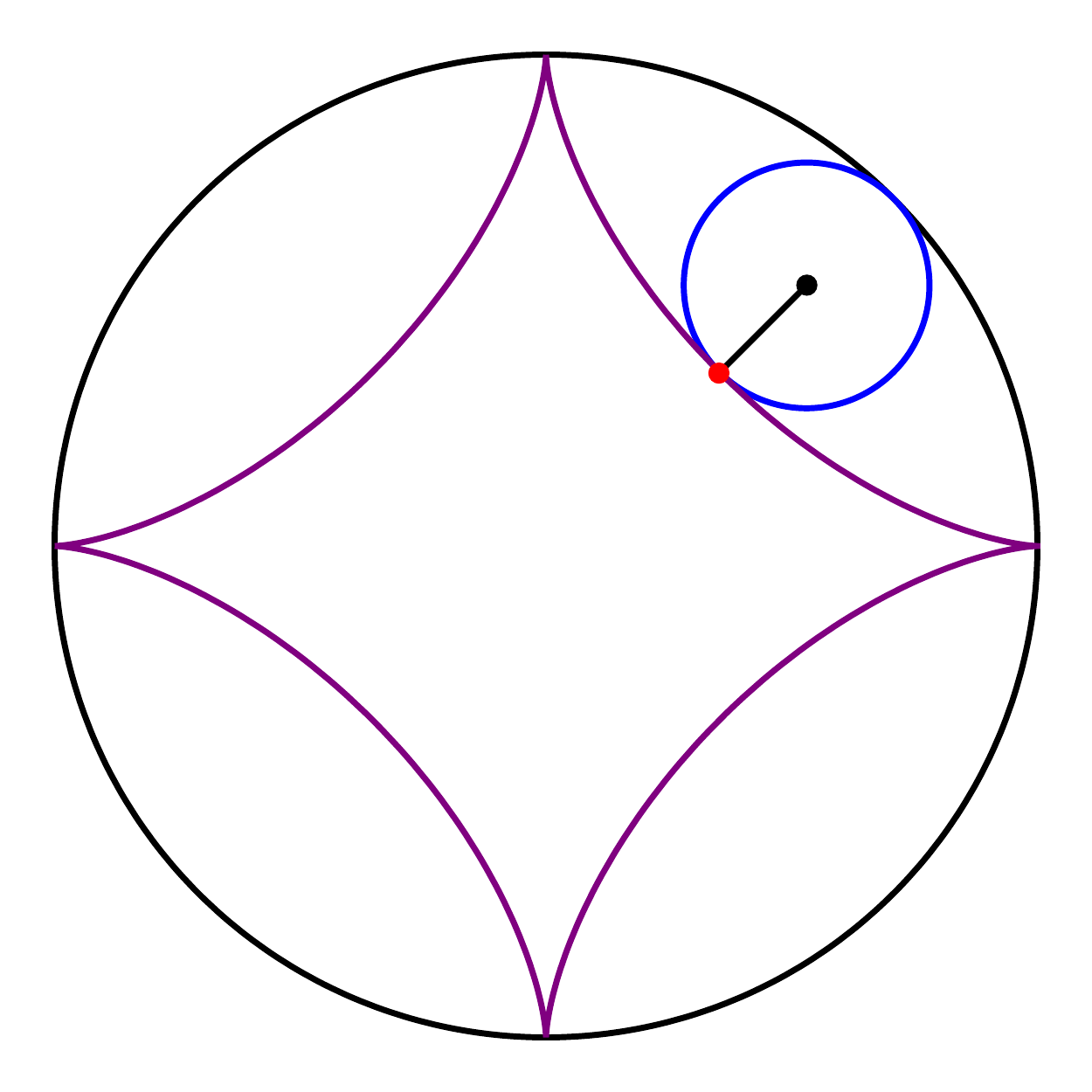}
		\caption{Astroid}
	\end{subfigure}
	\caption{Circles of radius 1 trace out hypocycloids as they roll within circles of radii (from left to right) 2, 3, and 4.}
	\label{fig:roll}
\end{figure}

A novel and particularly accessible behavior occurs when $d=r$ is an odd prime. The reader might recall that a \emph{hypocycloid} is a planar curve obtained by tracing a fixed point on a circle of integral radius as it ``rolls'' within a larger circle of integral radius. Figure~\ref{fig:roll} illustrates this construction. We are interested in the hypocycloid that is centered at the origin and has $r$ cusps, one of which is at $r$. This curve is obtained by rolling a circle of radius $1$ within a circle of radius $d$; it has the parametrization $\theta \mapsto (r-1)e(\theta)+e((1-r)\theta)$. Let $H_r$ denote the compact region bounded by this curve.

\begin{Proposition}\label{prop:hypo}
If $r$ is an odd prime, then the image of $g_r$ is $H_r$.
\begin{proof}
Since
\[
	\Phi_r(x)=x^{r-1}+x^{r-2}+\cdots +x+1,
\]
we obtain the formula
\[
	g_r(z_1,z_2,\ldots,z_{r-1})=z_1+z_2+\cdots+z_{r-1}+\frac{1}{z_1z_2\cdots z_{r-1}}.
\]
The image of this map is the set of all traces of matrices in $SU(r)$, the group of $r\times r$ complex unitary matrices with determinant $1$. This set is none other than $H_r$~\cite[Theorem 3.2.3]{Cooper}. In particular, the image under $g_d$ of the diagonal $z_1=z_2=\cdots=z_{r-1}$ is the boundary of $H_r$.
\end{proof}
\end{Proposition}

To expand on the previous example, suppose again that $r$ is an odd prime and $b$ is a positive integer. We have
\[
	\Phi_{r^b}(x)=\sum_{j=0}^{r-1} x^{jr^{b-1}},
\]
whence
\begin{equation}
	g_{r^b}(z_1,z_2,\ldots,z_{r^b-r^{b-1}})
	=
	\sum_{j=1}^{r^b-r^{b-1}} z_j+\sum_{j=1}^{r^{b-1}}
	\prod_{\ell=0}^{r-2} z_{j+\ell r^{b-1}}^{-1}.
	\label{eq:primepower}
\end{equation}
If $r=3$ and $b=2$, for instance, then the map is given by
\begin{equation*}
	g_9(z_1,z_2,z_3,z_4,z_5,z_6)=
	z_1+z_4+\frac{1}{z_1z_4}
	+z_2+z_5+\frac{1}{z_2z_5}
	+z_3+z_6+\frac{1}{z_3z_6}.
	\label{eq:fillout}
\end{equation*}

\begin{figure}[ht]
	\begin{subfigure}[b]{0.484\textwidth}
		\includegraphics[width=\textwidth]{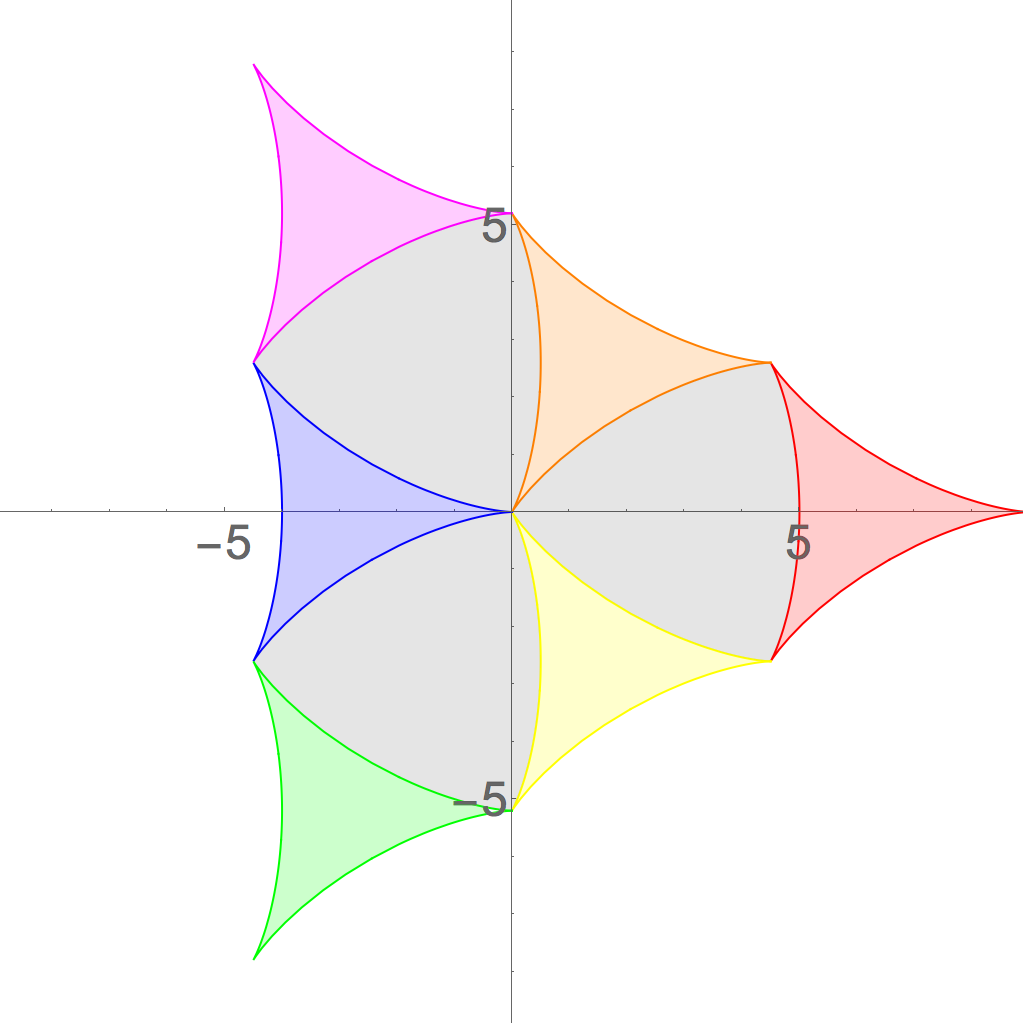}
		\caption{\scriptsize A geometric interpretation of $H_3+H_3+H_3$}
	\end{subfigure}
	\quad
	\begin{subfigure}[b]{0.484\textwidth}
		\includegraphics[width=\textwidth]{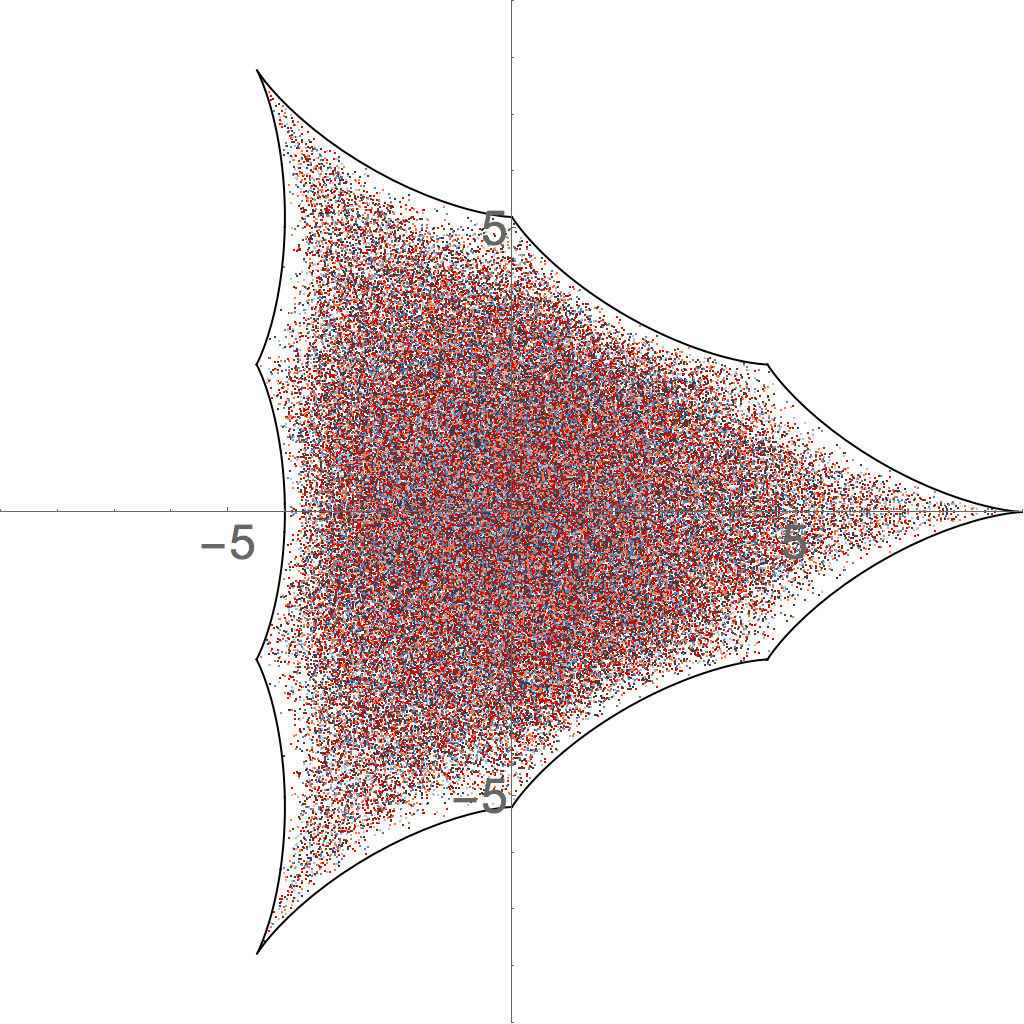}
		\caption{\scriptsize $n=1009^2$, $\omega=84669$, $d=3^2$}
	\end{subfigure}
	\caption{The supercharacter on the right fills out a Minkowski sum of filled deltoids.}
	\label{fig:mink}
\end{figure}

A definition will enable us to discuss the image in this situation. The \emph{Minkowski sum} of two nonempty subsets $S$ and $T$ of $\C$, denoted by $S+T$, is the set
\[
	S+T=\{s+t \in \C : s\in S\mbox{ and } t\in T\}.
\]
The Minkowski sum of an arbitrary finite collection is defined by induction. As a consequence of Proposition~\ref{prop:hypo}, we discover that the image of $g_9$ is none other than the Minkowski sum $H_3 + H_3 + H_3$, as illustrated in Figure~\ref{fig:mink}. A close look at~\eqref{eq:primepower} reveals a more general phenomenon.

\begin{Corollary}\label{cor:sumset}
If $r^b$ is a power of an odd prime, then the image of $g_{r^b}$ is the Minkowski sum
\[
	\sum_{j=1}^{r^{b-1}} H_r.
\]
\end{Corollary}

The Shapley--Folkman--Starr Theorem, familiar to mathematical economists, gives an explicit upper bound on the distance between points in a Minkowski sum and its convex hull \cite{Starr}. In the context of Corollary~\ref{cor:sumset}, we obtain the bound
\[
	\min\big\{ |w-z| : w\in \sum_{j=1}^{r^{b-1}} H_r\big\}
	\leq
	2\sqrt{2} r\sin\(\frac{\pi}{r}\),
\]
for any $z$ in the filled $r$-gon with vertices at $r^b e(\frac{j}{r})$ for $j=1,2,\ldots,r$. It follows easily that as $b\to\infty$, any point in the filled $r$-gon whose vertices are the $r$th roots of unity becomes arbitrarily close to points in the scaled Minkowski sum
\[
	\frac{1}{r^{b-1}}\sum_{j=1}^{r^{b-1}} H_r.
\]


\section{Concluding remarks}

There is work to be done toward understanding the images of cyclic supercharacters. If we are to stay the course of inquiry set by Section~\ref{sec:ex}, then a different approach is required; beyond the special cases discussed above, a general formula for the integers $c_{jk}$ in~\eqref{eq:cjk} appears unobtainable, since there is no known simple closed-form expression for the coefficients of an arbitrary cyclotomic polynomial $\Phi_d(x)$.

There is, however, a remedy. To minimize headache, suppose that $d=rs$ is a product of distinct odd primes, and that $\omega_q\mapsto (\gamma_r,\gamma_s)$ under the standard isomorphism $(\Z/d\Z)^\times \to (\Z/r\Z)^\times \times (\Z/s\Z)^\times$. Instead of wielding the elements $1,\omega_q,\ldots,\omega_q^{\phi(d)-1}$ as a $\Z$-basis for $\Z[\omega_q]$, we can use an analogous basis for $\Z[\gamma_r,\gamma_s]$. After some computation, we see that the image of the function $g_d$, formerly quite mysterious, is equal to the image of the function $h_d:\T^{\phi(d)}\to\C$ given by
\begin{equation}
	h\big((z_{ij})_{0\leq i<r-1,\,0\leq j<s-1}\big)=\sum_{i=0}^{r-2}\sum_{j=1}^{s-2} z_{ij} + \sum_{i=0}^{r-2}\prod_{j=0}^{s-2} \frac{1}{z_{ij}}+\sum_{j=0}^{s-2}\prod_{i=0}^{r-2} \frac{1}{z_{ij}} + \prod_{i=0}^{r-2}\prod_{j=0}^{s-2} z_{ij}.
	\label{eq:hfunc}
\end{equation}

This procedure, which amounts to a change of coordinates, can be used to obtain a closed formula for a Laurent polynomial map $h_d$ having the same image as $g_d$, for any integer $d$. This brings us one step closer to understanding the asymptotic behavior of the cyclic supercharacters in Section~\ref{sec:asym}. In practice, however, the functions $h_d$ are still difficult to analyze, despite being considerably easier to write down than the $g_d$. Even the simplest cases, described in~\eqref{eq:hfunc}, resist the accessible geometric description provided for the examples in Section~\ref{sec:ex}.

While certain graphical features of cyclic supercharacters with composite moduli have been explained in \cite{DukeGarciaLutz}, the mechanisms behind many of the striking patterns herein remain enigmatic. An important step toward deciphering the behavior of these supercharacters is to predict the layering constant $c$, discussed in Section~\ref{sec:cycsup}, given only a modulus $n$ and generator $\omega$. As is apparent, these layerings betray an underlying geometric structure that allows us to decompose the images of $\sigma_X$ into more manageable sets. We have been successful so far in finding appropriate values of $c$ ad hoc; however, a general theory is necessary to formalize our intuition.


\bibliographystyle{plainurl}
\bibliography{GHMFC2S}

\end{document}